\documentclass{amsart}
\usepackage{amsmath,amscd, amssymb,amsthm,amsfonts,graphicx}
\newcommand{\comment}[1]{}
\newcommand{\mathbold}[1]{{\mbox{{\boldmath${{#1}}$\unboldmath}}}}

\newcommand{\0}{\varnothing}
\newcommand{\defterm}[1]{{\it #1}}
\newcommand{\dju}{\sqcup}

\newcommand{\fld}{{\mathbb F}}
\newcommand{\gradient}{\nabla}
\newcommand{\half}{\frac{1}{2}}

\newcommand{\isom}{\cong}

\newcommand{\mdel}{\backslash}  
\newcommand{\normal}{\eta}
\newcommand{\ontopof}[2]{ \begin{matrix} #1 \\ \quad #2 \end{matrix} }
\newcommand{\ov}[1]{\overline{#1}}

\newcommand{\qbin}[2]{\left[ \begin{matrix} #1 \\ #2 \end{matrix} \right]_q}
\newcommand{\qbinsmall}[2]{\left[ \begin{smallmatrix} #1 \\ #2 \end{smallmatrix} \right]_q}
\newcommand{\sgn}{\varepsilon}

\newcommand{\sm}{\setminus}
\newcommand{\x}{\times}

\DeclareMathOperator{\Hom}{Hom}

\newcommand{\DLam}{{\mathcal L}^{d}}
\newcommand{\DRig}{{\mathcal R}^{d}}

\newcommand{\PhMap}{\varphi}           
\newcommand{\II}{{\mathcal I}}

\newcommand{\R}{{\mathcal R}}
\newcommand{\C}{{\mathcal C}}
\newcommand{\A}{{\mathcal A}}
\newcommand{\LL}{{\mathcal L}}
\newcommand{\HH}{{\mathcal H}}
\newcommand{\SSS}{{\mathcal S}}
\newcommand{\PP}{{\mathcal P}}

\newcommand{\Umat}[2]{U_{#1,#2}}

\newcommand{\Gr}{\mathbb{G}r}
\newcommand{\Zz}{\mathbb{Z}}
\newcommand{\Cc}{\mathbb{C}}
\newcommand{\Qq}{\mathbb{Q}}
\newcommand{\Rr}{\mathbb{R}}

\newcommand{\Pp}{\mathbb{P}}

\theoremstyle{plain}
\newtheorem{thm}{Theorem}[section]

\newtheorem{lem}[thm]{Lemma}
\newtheorem{cor}[thm]{Corollary}
\newtheorem{prop}[thm]{Proposition}

\theoremstyle{definition}
\newtheorem{defn}[thm]{Definition}
\newtheorem{ex}[thm]{Example}

\newtheorem{qn}[thm]{Question}
\newtheorem{problem}{Problem}

\theoremstyle{remark}
\newtheorem{rmk}[thm]{Remark}

\begin{document}

\title{Rigidity theory for matroids}
\author{Mike Develin, Jeremy L. Martin and Victor Reiner}
\address{Mike Develin, American Institute of Mathematics, 360 Portage
Ave., Palo Alto, CA 94306-2244, USA}
\email{develin@post.harvard.edu}
\address{Jeremy L. Martin, Department of Mathematics, University of Kansas,
Lawrence, KS 66045, USA}
\email{jmartin@math.ku.edu}
\address{Victor Reiner, School of Mathematics, University of Minnesota, Minneapolis, MN
55455, USA}
\email{reiner@math.umn.edu}
\date{\today}
\thanks{First author supported by the American Institute of
Mathematics.  Second author partially supported by an
NSF Postdoctoral Fellowship.  Third author partially supported by NSF grant
DMS--0245379}



\begin{abstract}
Combinatorial rigidity theory seeks to describe the rigidity or flexibility of
bar-joint frameworks in ${\mathbb R}^d$ in terms of the structure of the underlying
graph $G$.  The goal of this article is to broaden the foundations of
combinatorial rigidity theory by replacing $G$ with an arbitrary representable
matroid $M$.  The ideas of rigidity independence and parallel
independence, as well as Laman's and Recski's combinatorial characterizations
of 2-dimensional rigidity for graphs, can naturally be extended to this wider setting.
As we explain, many of these fundamental concepts really
depend only on the matroid associated with $G$ (or its Tutte polynomial), and
have little to do with the special nature of graphic matroids or the field ${\mathbb R}$.

Our main result is a ``nesting theorem'' relating the various
kinds of independence.
Immediate corollaries include generalizations of Laman's Theorem, as well as
the equality of 2-rigidity and 2-parallel independence.
A key tool in our study is the \emph{space of photos of $M$},
a natural algebraic variety whose irreducibility
is closely related to the notions of rigidity
independence and parallel independence.

%
%
The number of points on this variety, when working over a finite field,
turns out to be an interesting Tutte polynomial evaluation.
\end{abstract}

\maketitle

\tableofcontents

\section{Introduction: a brief tour through rigidity theory}

\emph{Combinatorial rigidity theory} is concerned with frameworks
built out of bars and joints in $\Rr^d$, representing the vertices $V$
and edges $E$ of an (undirected, finite) graph $G$.  (For comprehensive
treatments of the subject, see, e.g., \cite{GSS, Whiteley, Whiteley2}.)
The motivating problem is to
determine how the combinatorics of $G$ governs the rigidity
or flexibility of its frameworks.
Typically, one makes a generic choice of coordinates
  \begin{equation} \label{generic-coords}
  p=\{p_v:v \in V\} \subset \Rr^d
  \end{equation}
for the vertices of $G$, and considers
infinitesimal motions $\Delta p$ of the vertices.
The following two questions are pivotal:
\begin{enumerate}
  \item[(I.)] What is the dimension of the space of infinitesimal motions
    $\Delta p$ that preserve all the (squared) edge lengths $Q(p_u-p_v)$, for $\{u,v\}
    \in E$, where $Q(x)=\sum_{i=1}^d x_i^2$?
  \item[(II.)] What is the dimension of the space of infinitesimal motions
    $\Delta p$ that preserve all the edge directions $p_u-p_v$ regarded
    as slopes, that is, up to scaling?
\end{enumerate}

The answers to these questions are known to be determined by certain
linear dependence matroids represented over transcendental extensions of 
$\Rr$, as we now explain.

First, the $d$-dimensional {\it rigidity matroid} $\R^d(G)$ is the
matroid represented by the vectors
  \begin{equation}
  \label{graph-rigidity-vectors}
  \{ (e_u - e_v) \otimes (p_u-p_v) :\ \{u,v\} \in E\}
  \end{equation}
lying in $\Rr^{|V|} \otimes \Rr(p)^d$,
where $\Rr(p)$ is the extension of $\Rr$ by a collection of $d|V|$ transcendentals $p$,
thought of as the coordinates
of a generic embedding as in \eqref{generic-coords}.  The $|E| \times d|V|$
{\it rigidity matrix} $R^d(G)$ has as its rows the $|E|$ vectors in 
\eqref{graph-rigidity-vectors}.  Then the nullspace of $R^d(G)$ is
the space of infinitesimal motions of the vertices that preserve edge distances
(because $R^d(G)$ is $\frac{1}{2}$ times the Jacobian in the variables $p$
of the vector of squared edge lengths $Q(p_u-p_v)$; cf. 
Remark~\ref{Jacobian-interpretation} below).
Since row rank equals column rank, knowing the matroid $\R^d(G)$
represented by the rows of $R^d(G)$ answers question (I).

Second, the $d$-dimensional {\it parallel matroid} $\PP^d(G)$ is the
matroid represented by the vectors 
  \begin{equation}
  \label{graph-parallel-vectors}
  \{ (e_u - e_v) \otimes \normal_{u,v}^{(j)} :\ \{u,v\} \in E,\ j=1,2,\ldots,d-1\}
  \end{equation}
where for each edge $\{u,v\}\in E$, the vectors
$\normal_{u,v}^{(1)},\dots,\normal_{u,v}^{(d-1)}$
are generically chosen normals to $p_u-p_v$ in $\Rr^d$.
Again, we should consider the vectors in \eqref{graph-parallel-vectors}
as lying in $\Rr^{|V|} \otimes \Rr(p,\normal)^d$,
where $\Rr(p,\normal)$ is an extension of $\Rr$ by $d|V|$ transcendentals
$p$ and $(d-1)|E|$ transcendentals $\normal$.  In analogy to the preceding
paragraph, the $|E| \times d|V|$
\emph{parallel matrix} $P^d(G)$ has as its rows the $|E|$ vectors in
\eqref{graph-parallel-vectors}, and its nullspace
is the space of infinitesimal motions of the vertices that preserve
all edge directions.
Consequently, the matroid $\PP^d(G)$ represented by the rows of $P^d(G)$
provides the answer to question (II).

Some features of the theory are as follows:
\begin{enumerate}
\item[$\bullet$]  For $d=1$, the rigidity matroid coincides with the usual
graphic matroid for $G$ (while the parallel matroid is a trivial object).
\item[$\bullet$]  For $d=2$, the rigidity and parallel matroids coincide
\cite[Corollary~4.1.3]{Whiteley}.  Furthermore, this matroid 
$\R^2(G)=\PP^2(G)$ has many equivalent combinatorial reformulations, of which
the best known is {\it Laman's condition} \cite{Laman}:
$A \subseteq E$ is 2-rigidity-independent
if and only if for every subset $A' \subseteq A$
 \begin{equation}
  \label{Laman's-condition}
  \begin{aligned}
    2|V(A')| - 3 &\geq |A'|,\quad \text{ or equivalently} \\
    2\left(|V(A')| - 1\right) &> |A'|
  \end{aligned}
  \end{equation}
where $V(A')$ denotes the set of vertices incident to at least one edge in $A'$.
We will refer to this coincidence between $\R^2(G), \PP^2(G)$ and
the matroid defined by Laman's condition as the {\it planar trinity}.

\item[$\bullet$] For $d>2$, the parallel matroid has a simple 
combinatorial characterization that generalizes Laman's condition,
while an analogous description for the rigidity matroid is not known.
\end{enumerate}

\section{Main definitions: from graphs to matroids}

The purpose of this article is to broaden the scope of rigidity theory
by replacing the graph $G$ with a more general object: a matroid $M$
represented over an arbitrary field $\fld$.  As we shall see, the
notions of rigidity and parallel independence, as well as Laman's
combinatorial characterization, admit natural generalizations
to the setting of matroids.

In the process, we will see that many of the main results of
do not depend on the special properties of graphs (or graphic matroids),
nor on the field $\Rr$, but indeed remain valid for any matroid $M$
represented as above.  In the 
process, we are led naturally to an algebraic variety called the
\emph{space of $k$-plane-marked $d$-photos} of $M$.  Just as a
bar-joint framework may be regarded as an embedding of a graph in $\Rr^d$,
a photo of $M$ is a ``model'' of $M$ in $\fld^d$.

Whether or not the photo space is irreducible plays a key role in
characterizing the matroid analogues of rigidity independence
and parallel independence.  In turn, the question of irreducibility
can be answered combinatorially, using some elementary algebraic 
geometry and the classic matroid partitioning result of
Edmonds \cite{Edmonds}. We note in addition that when the field 
$\fld$ is finite, the number of photos of $M$ is counted by an
evaluation of the Tutte polynomial using $q$-binomial coefficients.

In order to summarize our results, we define the
main protagonists here.  Recall that a simplicial complex on vertex set
$E$ is a collection $\II$ of subsets of $E$ satisfying 
the following hereditary condition:
    if $I \in \II$ and $I' \subseteq I$, then $I' \in \II$.  
The independent sets of a matroid always form a simplicial complex. 
From here on we will make free use of standard terminology and notions
from matroid theory; background and definitions may be found in standard
texts such as \cite{Aigner, Oxley, Welsh}.

\begin{defn}[\textbf{$m$-Laman independence, $m$-Laman complex}]
\label{laman-matroid}
Let $E$ be a set of cardinality $n$, and let $M$ be a (not necessarily
representable) matroid on ground set $E$, with rank function $r$.  For
$m$ a real number in the open interval $(1,\infty)_\Rr$, say that $A
\subseteq E$ is {\it $m$-Laman independent} if 
  \begin{equation}
  \label{d-Laman-independent}
    m \cdot r(A') > |A'| \quad\text{ for all nonempty subsets } A' \subseteq A.
  \end{equation}
The \defterm{$m$-Laman complex} $\LL^m(M)$ is defined
as the abstract simplicial complex
of all $m$-Laman independent subsets of $E$.
\end{defn}

We will prove combinatorially that
\begin{enumerate}
  \item[$\bullet$] If $m$ is a positive integer, then $\LL^m(M)$ is the
    collection of independent sets of a matroid (Theorem~\ref{Laman-matroidal}),
    but this is not true in general for other values of $m$.
  \item[$\bullet$] If $m$ is a positive integer, then $\LL^m(M)$ has several
    other combinatorial characterizations (Theorem~\ref{Edmonds-Laman-equiv}), 
    including a generalization of Recski's Theorem.
  \item[$\bullet$] If $m=2$ and $M$ is representable, then
    the matroid $\LL^2(M)$ coincides with the $2$-dimensional
    rigidity and parallel matroids, defined below (Corollary~\ref{planar-case}).
\end{enumerate}

Throughout the rest of the introduction, let $M$ be a represented matroid;
that is, a matroid equipped with a representation over some field $\fld$ by vectors
$E=\{v_1,\ldots,v_n\}\subset\fld^r$.  It is worth emphasizing that we are
\emph{not} regarding $M$ as an abstract matroid; that is, the vectors
$\{v_1,\dots,v_n\}$ are part of the data of $M$.  For notational convenience,
we identify the ground set $E$ with the numbers $[n]:=\{1,2,\ldots,n\}$.
Denote by $\Gr(k,\fld^d)$ the Grassmannian of $k$-planes in $\fld^d$, regarded
as a projective variety over $\fld$ via the usual Pl\"ucker embedding into
$\Pp^{\binom{d}{k}-1}$.

When $m>1$ is a rational number, the Laman complex 
$\LL^m(M)$ is closely related to an algebraic variety that we now define.

\begin{defn}[\textbf{photo space, $(k,d)$-slope independence, $(k,d)$-slope complex}]
\label{photo-space}
Let $M$ be a matroid equipped with representation $\{v_1,\dots,v_n\}$ as above.
The corresponding \defterm{space of $k$-plane-marked $d$-photos} (or just \defterm{$(k,d)$-photos})
is the algebraic set
  \begin{equation}
  X_{k,d}(M):= \{ (\PhMap,W_1,\dots,W_n) \in \Hom_\fld(\fld^r,\fld^d) \times \Gr(k,\fld^d)^n:~
    \PhMap(v_i) \in W_i \text{ for } i=1,\ldots,n\}.
  \end{equation}

The photo space of a matroid is analogous to the picture space of a
graph, as defined in \cite{JLM1,JLM2}.
One may think of the map $\PhMap \in \Hom_\fld(\fld^r,\fld^d)$ 
as projecting the vectors $\{ v_i \}$ into a space $\fld^d$ of dimension
possibly less than $r$, like a camera taking a photo of the $\{ v_i\}$ on 
photographic paper that looks like 
$\fld^d$.  The $k$-plane $W_i$ in $\fld^{d}$ is thought of as a ``marking'' of the image vector
$\PhMap(v_i)$ in the photo, so that $W_i$ is constrained to contain $\phi(v_i)$.
Of course, whenever $\PhMap(v_i)=0$ (perhaps the camera $\PhMap$ caught $v_i$ at
a bad angle), this $k$-plane $W_i$ is unconstrained.  
The idea of $(k,d)$-slope independence is to measure how independently these marking $k$-planes
can vary while obeying these constraints, when none of the $v_i$ are annihilated by $\PhMap$.
The linear dependences among the $v_i$ force linear dependences among their 
image vectors $\PhMap(v_i)$, and hence algebraic constraints among the subspaces $W_i$. 

Define a Zariski open subset of $X_{k,d}(M)$ (called the
\defterm{non-annihilating cellule};
see Definition~\ref{cellule-definition} below) by
  $$
  X^{\0}_{k,d}(M):= \{ (\PhMap,W_1,\dots,W_n) \in X_{k,d}(M):~
    \PhMap(v_i) \neq 0 \text{ for } i=1,2,\ldots,n \}
  $$
and consider its image under the projection map
  \begin{equation} \label{projection-map}
    \Hom(\fld^r,\fld^d) \times \Gr(k,\fld^d)^n \ 
    \overset{\pi}{\longrightarrow} \ \Gr(k,\fld^d)^n.
  \end{equation}
This image measures the constraints on the $W_i$ when none of the $v_i$
are mapped to zero;
specifically, we define $M$ to be \defterm{$(k,d)$-slope independent}
if $\pi X^{\0}_{k,d}(M)$
is Zariski dense in $\Gr(k,\fld^d)^n$.  The \defterm{$(k,d)$-slope complex}
is defined as
  \begin{equation}
  \label{slope-complex-defn}
  \SSS^{k,d}(M):= \{ A \subseteq E: M|_A \text{ is }(k,d)\text{ -slope independent}\}.
  \end{equation}
\end{defn}

A third notion of matroid rigidity generalizes the
$d$-dimensional rigidity matroid $\R^d(G)$ of a graph $G$.

\begin{defn}[\textbf{rigidity matroid, rigidity complex}]
\label{rigidity-matroid}
Let $M$ be a matroid equipped with representation $\{v_1,\dots,v_n\}$ as above,
and let $d$ be a positive integer.
The \defterm{$d$-dimensional (generic) rigidity matroid} 
is the matroid represented by the vectors
  \begin{equation}
  \label{general-rigidity-vectors}
  \{v_i \otimes \PhMap(v_i)\}_{i=1}^n \subset \fld^r \otimes_\fld \fld(\PhMap)^d.
  \end{equation}
where $\fld(\PhMap)$ is the field extension of $\fld$ by $dr$ transcendentals giving
the entries of the matrix $\PhMap: \fld^r \rightarrow \fld(\PhMap)^d$.
The \defterm{$d$-rigidity complex} $\R^d(M)$ is the complex of independent sets of the
$d$-dimensional rigidity matroid, and the
\defterm{$d$-rigidity matrix} $R^d(M)$ is the $n \times dr$ matrix whose rows 
are given by the vectors $\eqref{general-rigidity-vectors}$.
\end{defn}

To interpret this construction, consider
the pseudo-distance quadratic form $Q(x):=\sum_{i=1}^d x_i^2$ on $\fld(\PhMap)^d$.
Provided that the field $\fld$ has characteristic $\neq 2$, one can interpret 
the nullspace of $R^d(M)$ as the space of infinitesimal changes 
of $\PhMap$ that
preserve the values $Q(\PhMap(v_i))$ for all $i=1,2,\ldots,n$. (See 
Proposition~\ref{rigidity-parallel-interpretations}(ii).)

\begin{defn}[\textbf{hyperplane-marking matroid}]
\label{hyperplane-marking-matroid}
Let $M$ be a matroid represented by $v_1,\ldots,v_n \in \fld^r$
as above.  Its {\it ($d$-dimensional, generic) hyperplane-marking matroid} 
is the matroid represented over $\fld(\PhMap, \normal)$
by the vectors
  $$
  \{v_i \otimes \normal_i\}_{i=1}^n \subset \fld^r \otimes_\fld \fld(\PhMap,\normal)^d
  $$
where $\fld(\PhMap,\normal)$ is the extension of $\fld$ by $dr$
transcendentals $\PhMap_{ij}$
(the entries of the matrix $\PhMap$) and $(d-1)n$ more transcendentals $\normal_{ij}$
(the coordinates of the vectors $\normal_i$ normal to $\PhMap(v_i)$). 
The complex $\HH^d(M)$ is defined to be the complex of independent sets of 
this matroid.
\end{defn}

To interpret the notion of rigidity independence modeled by $\HH^d(M)$, one
should regard lack of rigidity as the ability to deform $\PhMap$ so that
the images $\PhMap(v_i)$ of the ground set elements vary, but membership 
in their orthogonal complement hyperplanes is preserved.
The most important instance of the hyperplane-marking matroid
uses the \defterm{$(d-1)$-parallel extension} of $M$,
the matroid $(d-1)M$ whose ground set consists of $d-1$
parallel copies of each element of $E$.  The
\defterm{($d$-dimensional, generic) parallel matroid} is defined as
  $$\PP^d(M) := \HH^d((d-1)M),$$
and the \defterm{$d$-parallel matrix} $P^d(M)$ is defined as the
$n \times dr$ matrix whose rows represent $\HH^d((d-1)M)$.  Its nullspace
consists of the infinitesimal changes $\Delta \PhMap$ in the matrix $\PhMap$ which
preserve the slopes of all the direction vectors $\PhMap(v_i)$
(see Proposition~\ref{rigidity-parallel-interpretations}~(i)).

These definitions generalize the ordinary definitions from the rigidity 
theory of graphs.  Strikingly, the geometric constraints 
on the photo space can be categorized combinatorially: the identity
  $$
  \SSS^{k,d}(M) = \LL^{\frac{d}{d-k}}(M),
  $$
(Corollary~\ref{Laman-slope-relation})
provides a geometric interpretation of $\LL^m(M)$ for rational $m$.

The slope complex $\SSS^{k,d}(M)$ is closely related to the rigidity and 
parallel matroids.  The precise relationship is given by the \emph{Nesting Theorem}
(Theorem~\ref{Laman-rigidity-relation}):
  $$
  \SSS^{1,d}(M) \subseteq \R^d(M) \subseteq \LL^d(M) = \HH^d(M) = \SSS^{d-1,d}(M)
  $$
for all integers $d\geq 2$.  In particular, when $d=2$,
  \begin{equation} \label{plane-case}
  \HH^2(M) = \SSS^{1,2}(M) = \R^2(M) = \LL^2(M).
  \end{equation}
Thus matroid rigidity theory leads to a conceptual proof of 
the planar trinity (the second and third inequalities in \eqref{plane-case}).

For $d \geq 3$, the $d$-rigidity matroid $\R^d(M)$ is the hardest of these
objects to understand (as it is for graphic matroids).
One fundamental question is whether $\R^d(M)$ depends on the choice of representation
of $M$.  It is invariant for $d=2$ (by the Nesting Theorem) and
up to projective equivalence of representations
(Proposition~\ref{projective-invariance}), but the problem remains open for
the general case.  We also study the behavior of
the $d$-rigidity matroid as $d\to\infty$, and show (Proposition~\ref{stabilization})
that $R^d(M)$ stabilizes when $d\geq r(M)$.

\section{Laman independence}
\label{Laman-section}

The central result of this section, Theorem~\ref{Laman-matroidal}, states that
the generalized Laman's condition \eqref{d-Laman-independent} always gives
a matroid when $m$ is an integer.  The proof is completely combinatorial;
that is, it is a statement about abstract matroids, not represented matroids.
In addition, we describe some useful equivalent characterizations
of $d$-Laman independence: one uses the Tutte polynomial, another is reminiscent
of Recski's Theorem, and another is related to Edmonds' theorem
on decomposing a matroid into independent sets.

\subsection{When is the Laman complex matroidal?}

\begin{thm}
\label{Laman-matroidal}
\begin{enumerate}
\item[(i)] 
    Let $d$ be a positive integer and let $M$ be any matroid.
    Then the simplicial complex $\DLam(M)$ is a matroid complex.
\item[(ii)]
    Let $m \in (1,\infty)_\Rr$ be a real number which is {\it not} an integer.
    Then there exists a represented
    matroid $M$ for which $\LL^m(M)$ is not a matroid complex.
\end{enumerate}
\end{thm}

\begin{proof}
For the first assertion, it is most convenient to use the characterization of matroids
by circuit axioms \cite[eq.~6.13,~p.~264]{Aigner}.
Define $\C$ to be the collection of those subsets of $E$ which are
minimal among nonmembers of $\DLam(M)$.
We wish to show that $\C$ satisfies the axioms for the circuits of a matroid.
Since $\DLam(M)$ is a simplicial complex, we only need check the {\it circuit exchange
axiom}:

\begin{center}
if $C, C' \in \C$ with $C \neq C'$, and $e \in C \cap C'$,
then there exists $C'' \in \C$ such that $C'' \subseteq (C \cup C') \sm \{e\}$.
\end{center}

Since $C, C'$ are minimal among the sets not satisfying the hereditary property
\eqref{d-Laman-independent}, we claim that
  $$
  \begin{aligned}
  |C|  &= d \cdot r(C), \\
  |C'| &= d \cdot r(C'), \\
  \end{aligned}
  $$
where $r$ is the rank function of $M$.  To see this claim, note
that $|C| \geq d \cdot r(C)$, and if this inequality were strict, 
then 
$$
|C-e|\geq d\cdot r(C)\geq d\cdot r(C-e)
$$
for any $e\in C$, contradicting the statement that $C$ is a \emph{minimal}
set not satisfying \eqref{d-Laman-independent}.  
Note also that $C \cap C'$ is a proper subset of each of $C,C'$ and hence
$$
  |C \cap C'| < d \cdot r(C \cap C').
$$
Since $d$ is an integer, the last condition may be rewritten as 
  $$
  |C \cap C'|+1 \leq d \cdot r(C \cap C').
  $$
The rank submodular inequality $r(C \cup C') \leq r(C) + r(C') - r(C \cap C')$
then implies
  \begin{eqnarray*}
    d \cdot r((C \cup C') \sm \{e\}) &\leq& d \cdot r(C \cup C') \\
      &\leq& d \cdot r(C) + d \cdot r(C') - d \cdot r(C \cap C') \\
      &\leq& |C| + |C'| - |C \cap C'| - 1\\
      &=& |(C \cup C') \sm \{e\}|.
  \end{eqnarray*}
So $(C \cup C') \sm \{e\}$ is not in $\DLam(M)$, hence contains some element of
$\C$.  This establishes (i).
\bigskip

We now prove (ii).  Suppose that $m \in (1,\infty)_\Rr$ is {\it not} an integer,
and let $c:=\lfloor m \rfloor$ (the greatest integer $\leq m$).
Choose positive integers $a,b$
satisfying the inequalities
\eqref{range-for-m} in Lemma~\ref{interval-exhaustion} below.
We will explicitly construct a represented matroid $M_{a,b,c}$
such that $\LL^m(M_{a,b,c})$ is not a matroid complex.

Let $\fld$ be a sufficiently large (for example, infinite) field,
let $V$ be a $(2b-1)$-dimensional vector space over $\fld$, 
and let $V_1,V_2$ be two $b$-dimensional subspaces of $V$
whose intersection $V_1 \cap V_2 = \ell$ is a line.  
Let $X=\{x_1,\ldots,x_c\}$ be a set of $c$ nonzero vectors on $\ell$.
For $i=1,2$, choose a set $Y_i \subseteq V_i$ of
cardinality $a-c$ generically (this is always possible if $\fld$ is
sufficiently large).  
Note in particular that no member of $Y_1\cup Y_2$ lies on $\ell$.

Let $M_{a,b,c}$ be
the matroid represented over $\fld$ by $E=X \cup Y_1 \cup Y_2$, and
denote by $\C$ the set of subsets of $E$ that are minimal among
nonmembers of $\LL^m(M_{a,b,c})$.  We claim that $\C$ does not satisfy the
circuit exchange axiom.  To see this,
let $C_i = X \cup Y_i$ for $i=1,2$ and observe that 
  $$
  m\cdot r(C_i) = mb \leq a = |C_i|,
  $$
so $C_i \not\in \LL^m(M_{a,b,c})$.  In fact, we claim that $C_i \in \C$.
Indeed, let $I$ be any nonempty proper subset of $C_i$ and let
$J = I \cap Y_i$.  Since $r(X)=1$,
and by the generic choice of $Y_1$ and $Y_2$, we have
  \begin{eqnarray*}
    r(J) &=& \min(|J|,b),\\
    r(I) &=& \min(|J|+1,b),\\
    m\cdot r(I) &=& \min(m|J|+m,mb).
  \end{eqnarray*}
Now Lemma \ref{interval-exhaustion} implies that $mb \geq a = |C_i| > |I|$.
Since $m$ is not an integer, we have also
  $$m|J|+m > |J|+c ~=~ |J|+|X| \geq |I|.$$
In all cases $m \cdot r(I) > |I|$.  It follows that $C_i \in \C$.

Now, let $x_i \in X$, and let $I = (C_1 \cup C_2) \sm \{x_i\} = E \sm \{x_i\}$.
Then every nonempty subset $I' \subseteq I$ satisfies
\eqref{d-Laman-independent}.  (We omit the routine but tedious calculation, which
involves eight cases, depending on how $I'$ meets each of $X$, $Y_1$ and $Y_2$.)
That is, $I$ is $m$-Laman-independent, hence contains no element of $\C$.
Therefore $\C$ fails the circuit exchange axiom, and we are done.
\end{proof}

The following numerical lemma was used in the proof of Theorem~\ref{Laman-matroidal}.

\begin{lem}
\label{interval-exhaustion}
Let $m \in (1,\infty)_\Rr$ be a real number which is not an integer, and
let $c:=\lfloor m \rfloor$.  Then there exist positive integers $a,b$ such that
  \begin{equation}
  \label{range-for-m}
  \frac{a-1}{b} < \frac{2a-c-1}{2b-1} < m \leq \frac{a}{b}.
  \end{equation}
\end{lem}

\begin{proof}
First, note that the third inequality implies the first one.
Indeed, if $m \leq a/b$, then 
$$
b+a \geq 1+a \geq 1+bm > 1+bc,
$$
which implies in turn that $2ab-a-2b+1 < 2ab-bc-b.$  Factoring this gives
$(2b-1)(a-1) < b(2a-c-1)$, or $\frac{a-1}{b} < \frac{2a-c-1}{2b-1}$ as desired.

We therefore concentrate on the second and third inequalities in \eqref{range-for-m}.
Subtracting $c$ from each expression in \eqref{range-for-m}
and substituting $a=bc+r$ yields
  \begin{equation} \label{br-constraint}
  \frac{2r-1}{2b-1} < m-c \le \frac{r}{b} = \frac{2r}{2b}.
  \end{equation}
Therefore, it will suffice to find
a pair $b,r$ of positive integers satisfying \eqref{br-constraint}.

Note that $m-c$ is the fractional part of $m$; since $m$ is not an integer,
we have $m-c \in (0,1)_{\Rr}$.  Therefore, it will suffice to show that $(0,1)$ is
the union of intervals of the form $(\frac{2r-1}{2b-1},\frac{2r}{2b}]$
for positive integers $b, r$.  Indeed,
  \begin{align*}
  (0,1) &= \bigcup_{m \geq 0} \left(\frac{m}{m+1},\ \frac{m+1}{m+2}\right]\\
        &= \left(\frac{0}{1},\ \frac{1}{2}\right] 
            \cup \left(\frac{1}{2},\ \frac{2}{3}\right] 
             \cup \left(\frac{2}{3},\ \frac{3}{4}\right]
              \cup \cdots \\
\intertext{and}
  \left(\frac{m}{m+1}, \frac{m+1}{m+2}\right] 
        &= \bigcup_{s \geq 1} \left(\frac{2^s m+1}{2^s(m+1)+1},\ 
              \frac{2^s m +2}{2^s(m+1)+2}\right] \\
        &= \left(\frac{2m+1}{2m+3},\ \frac{2m+2}{2m+4}\right] 
            \cup \left(\frac{4m+1}{4m+5},\ \frac{4m+2}{4m+6}\right]
             \cup \left(\frac{8m+1}{8m+9},\ \frac{8m+2}{8m+10}\right] 
              \cup \cdots
  \end{align*}
establishing \eqref{br-constraint}, as desired.
\end{proof}

\subsection{Equivalent characterizations}

One of the equivalent phrasings of $m$-Laman independence involves the
\defterm{Tutte polynomial} $T_M(x,y)$ of $M$, a fundamental isomorphism
invariant of the matroid $M$.  For background on the Tutte polynomial,
see the excellent survey article by Brylawski and Oxley \cite{BryOx}.

Given a subset $A$ of the ground set
$E$, denote by $\ov{A}$ the {\it matroid closure} or {\it span} of $A$.
If $A=\ov{A}$, then $A$ is called a \defterm{flat} of $M$.

\begin{prop}
\label{flat-condition}
Let $M$ be a matroid on ground set $E$ with rank function $r$,
and fix $m \in (1,\infty)_\Rr$.

Then the following are equivalent:
  \begin{enumerate}
  \item[(i)] $E$ is $m$-Laman independent, that is, $\LL^m(M)=2^E$ (the power set of $E$). 
  \item[(ii)] $m \cdot r(\ov{A}) > |\ov{A}|$ for every nonempty subset $A \subseteq E$.
      (Equivalently, $m \cdot r(F) > |F|$ for every flat $F$ of $M$.)
  \item[(iii)] The Tutte polynomial specialization $T_M(q^{m-1},q)$ is monic
    of degree $(m-1)r(M)$.
  \end{enumerate}
\end{prop}
\noindent
Note that in (iii) we must allow (non-integral) real number exponents for a ``polynomial''
in $q$, but the notions of ``degree'' and ``monic'' for such polynomials should still
be clear.  The connection between the Tutte polynomial
and rigidity of graphs was observed by the second author in \cite[\S 6]{JLM2}.

\begin{proof}
The equivalence of (i) and (ii) is clear from
the definition of $m$-Laman independence since $r(\ov{A})=r(A)$ and $|\ov{A}| \geq |A|$
for any $A \subseteq E$.

For the equivalence of (i) and (iii)
we use Whitney's {\it corank-nullity} formula \cite[eq.~6.13]{BryOx}
for the Tutte polynomial:
  $$
  T_M(x,y) = \sum_{A \subseteq E} (x-1)^{r(M)-r(I)} (y-1)^{|I|-r(I)}
  $$
Substituting $x=q^{m-1}$ and $y=q$ gives an expression for
$T_M(q^{m-1},q)$ as a sum of terms indexed by subsets $A \subseteq E$,
each of which is a monic polynomial in $q$ of degree
  $$
  (m-1)r(M) - m \cdot r(A) + |A|.
  $$
Thus $T_M(q^{m-1},q)$ will have degree at most $(m-1)r(M)$ if and only if
$m \cdot r(A) \geq |A|$  for all subsets $A \subseteq E$.
Furthermore, since the term indexed by $A=\0$ is monic
of degree $(m-1)r(M)$, the whole polynomial $T_M(q^{m-1},q)$ will be monic
of degree $(m-1)r(M)$ if and only if $m \cdot r(A) > |A|$ for every
nonempty subset $A$, that is, if and only if $E$ is $m$-Laman independent.
\end{proof}

Suppose that $m=d$ is a positive integer, so that $\DLam(M)$ is a matroid 
complex. Here $d$-Laman independence has two more equivalent 
formulations, one of which extends a classical result in the rigidity 
theory of graphs.

\vskip 0.1in
\noindent {\bf Recski's Theorem} \cite{Recski}.
\textit{Let $G = (V,E)$ be a graph, and let $E'$ be a spanning set of edges of 
size $2|V|-3$. Then $E'$ is a 2-rigidity basis if and only if for any 
$e\in E'$, we can partition the multiset $E'\cup \{e\}$ (that is, adding
an extra copy of $e$ to $E'$) into two disjoint spanning trees of $G$.}
\vskip 0.1in

This notion can be naturally extended to arbitrary matroids and dimensions.

\begin{defn}\label{Recski-indep}
Let $M$ be a matroid on $E$. We say that $E$ is \defterm{$d$-Recski 
independent} if for any element $e\in E$, the multiset $E\cup \{e\}$ can 
be 
partitioned into $d$ disjoint independent sets for $M$.
\end{defn}

We wish to show that this purely matroidal condition is equivalent
to the purely matroidal condition of $d$-Laman independence. 
To prove this, we use a powerful classic 
result of Edmonds.

\vskip 0.1in
\noindent {\bf Edmonds' Decomposition Theorem}
\cite[Theorem 1]{Edmonds}.
\textit{Let $M$ be a matroid of rank $r$ on ground set $E$.
Then $E$ has a decomposition $E = I_1 \sqcup \cdots \sqcup I_d$ into disjoint
independent sets $I_j$ for each $j=1,\dots,d$ if and only if
$d \cdot r(A) \geq |A|$ for every subset $A \subseteq E$.}

\begin{defn} 
\label{Edmonds-decomp}
Let $M$ be a matroid on $E$.  A \defterm{$d$-Edmonds decomposition} of $M$ is a
family of independent sets $I_1,\dots,I_d$ whose disjoint union is $E$, with
the following property: given subsets $I'_1 \subseteq I_1$, \dots, $I'_d \subseteq I_d$
with not all $I'_i$ empty, then it is not the case that $\ov{I'_1} = \ov{I'_2}
= \dots = \ov{I'_d}$.
\end{defn}

\begin{thm}
\label{Edmonds-Laman-equiv}
Let $M$ be a matroid on ground set $E$, and let $d$ be a positive integer. 
Then the following are equivalent:

(i) $E$ has a $d$-Edmonds decomposition;

(ii) $E$ is $d$-Laman independent;

(iii) $E$ is $d$-Recski independent.

\end{thm}

\begin{proof}
{\sf (ii)~$\Rightarrow$~(i):}
Suppose that $E$ is $d$-Laman independent.  By Edmonds'
Theorem, the set $E$ can be partitioned into disjoint independent sets 
$I_1,\dots,I_d$.  We claim that every such family is a $d$-Edmonds 
decomposition.  Indeed, suppose that $I'_1 \subseteq I_1$, \dots, $I'_d 
\subseteq I_d$ all have the same span, with not all $I'_j$ empty.  Since the 
$I_j$ are independent, the $I'_j$ all have the same cardinality, say $s$.  
Then $A := I'_1 \dju \dots \dju I'_d$ is nonempty and has rank $s$ and 
cardinality $ds$, which violates the $d$-Laman independence of $E$.

\smallskip\noindent
{\sf (i)~$\Rightarrow$~(ii):}
Let $I_1,\dots,I_d$ be a $d$-Edmonds decomposition of $M$.  Let
$A \subseteq E$ be nonempty, and $A_j := A \cap I_j$.  Then
  $$
  |A| = \sum_{j=1}^d |A_j| = \sum_{j=1}^d r(A_j) \leq \sum_{i=1}^d r(A)
  = d \cdot r(A).
  $$
However, equality cannot hold:  it would force $r(A_j) = r(A)$ for each $j$, 
so that each $A_j$ has the same span as $A$, violating 
the definition of a $d$-Edmonds decomposition.  Hence $|A| < d \cdot r(A)$ as
desired.

\smallskip\noindent
{\sf (ii)~$\Rightarrow$~(iii):}
Suppose that $E$ is $d$-Laman independent.
Consider the matroid $M'$ given by cloning any $e\in 
E$ as in the definition of $d$-Recski independence, so that the ground set 
of $M'$ is $E' = E\cup \{e\}$. We claim that
$|A'|\le d\cdot r(A')$ for each $A'\subseteq E'$. Indeed, either
$A'\subseteq E$, when $|A'| < d\cdot r(A')$, or else $A' = 
A\cup \{e\}$ with $A\subseteq E$, when $|A'| = |A| + 1 < 
d\cdot r(A) + 1$, so $|A'| \le d\cdot r(A) \le d\cdot r(A')$.  By
Edmonds' Theorem, $E'$ can be partitioned into $d$ 
disjoint independent subsets.  It follows that $M$ is $d$-Recski independent.

\smallskip\noindent
{\sf (iii)~$\Rightarrow$~(ii):}
Suppose that $E$ is not $d$-Laman independent, i.e., it has a subset
$A$ with $|A| \geq d\cdot r(A)$. Let $a\in A$. The set
$A\cup \{a\}\subseteq E\cup \{a\}$ has rank $r(A)$ and cardinality 
$|A|+1$, so $|A\cup\{a\}| > d\cdot r(|A\cup\{a\}|)$. 
By Edmonds' Theorem, $E\cup\{a\}$ cannot be partitioned into $d$
independent sets.  Hence $E$ is not $d$-Recski independent.
\end{proof}

\subsection{Digression on polymatroids}

As we have seen in Theorem~\ref{Laman-matroidal}~(ii),
when $m$ is not an integer, the Laman complex $\LL^m(M)$
need not form the collection of independent sets of a matroid.
However, $\LL^m(M)$ is related to a more general (and less
well-known) object called a \defterm{polymatroid}, as we
now explain.  (The results of this section will not
be necessary for the remainder of the paper.)

We review the definition of a polymatroid,
using its characterizations by monotone submodular rank functions
(see \cite[Chapter 18]{Welsh}).

\begin{defn}
Fix the ground set $E=[n]$.  A function $\rho: 2^E \rightarrow \Rr_{\geq 0}$ is
the \defterm{rank function of a polymatroid on $E$} if it is 
\begin{itemize}
\item[--] {\it monotone}: $\rho(A) \leq \rho(B)$ whenever $A \subseteq B \subseteq E$;
\item[--] {\it submodular}: $\rho(A \cup B) + \rho(A \cap B) \leq \rho(A) + \rho(B)$
  for all $A,B \subseteq E$; and 
\item[--] {\it normalized}: $\rho(\0)=0$.
\end{itemize}

The {\it polymatroid} associated with $\rho$ is the convex polytope
  $$
  P_\rho:=\{ x \in \Rr_{\geq 0}^n:~ \sum_{a \in A} x_a \leq \rho(A)
    \text{ for all }A \subseteq E\},
  $$
also called the set of {\it independent vectors} of the polymatroid.
\end{defn}

Note that, for all $A \subseteq E$, the characteristic vector $\chi_A \in \Rr^n$ is
independent for $\rho$ if and only if $\rho(A) = |A|$.

Our goal is to show the following:

\begin{prop}
\label{Laman-polymatroidal}
For every loopless matroid $M$ on ground set $E=[n]$, and every real
number $m \in (1,\infty)_\Rr$,
there is a polymatroid rank function $\rho$ on $E$ with the following property:
$A\subseteq E$ is $m$-Laman independent if and only if its characteristic vector is
independent for $\rho$.
\end{prop}

The proof uses two standard lemmas.

\begin{lem} \cite[Lemma 6.15]{Aigner}
\label{Aigner-lemma}
If $f:2^E \rightarrow \Rr_{\geq 0}$ is monotone, submodular, and normalized, then so is
the function $r_f:2^E \rightarrow \Rr_{\geq 0}$ defined by
$$
r_f(A):=\min_{A' \subseteq A} \{ f(A') + |A \sm A'|  \}.
$$
\end{lem}

\begin{lem} \cite[Proposition A.3.1]{Whiteley}
\label{Whiteley-lemma}
For a monotone, submodular, normalized function $f:2^E \rightarrow 
\Rr_{\geq 0}$
with associate function $r_f$ as above, the following are equivalent:
\begin{enumerate}
\item[(a)] $|A'| \leq f(A') \text{ for all }A' \subseteq A.$
\item[(b)] $|A'| \leq r_f(A') \text{ for all }A' \subseteq A.$
\item[(c)] $r_f(A)=|A|.$
\end{enumerate}
\end{lem}

\begin{proof}[Proof of Proposition~\ref{Laman-polymatroidal}]
Let $\epsilon \in (0,\frac{1}{r(M)})_{\Rr}$, and define 
$f:2^E \rightarrow \Rr_{\geq 0}$ by
$$ 
f(A) = (m-\epsilon) r(A).
$$  
Note that $f$ is monotone, submodular, and normalized,
because the rank function $r$ of any loopless matroid has these properties.
By Lemma~\ref{Aigner-lemma}, the function $\rho:=r_f$ shares these properties,
hence also defines a polymatroid rank function on $E$.

Since $M$ is loopless, for all $\A \neq \0$,
one has $m r(A) > |A|$ if and only if $(m-\epsilon) r(A) \geq |A|.$
Consequently
  \begin{eqnarray*}
  A \in \LL^m(M) 
    &\iff& f(A') \geq |A'| \quad\text{ for all nonempty } A' \subseteq A \\
    &\iff& f(A') \geq |A'| \quad\text{ for all } A' \subseteq A \\
    &\iff& \rho(A) = r_f(A) = |A|.
  \end{eqnarray*}
Here the last equality uses Lemma~\ref{Whiteley-lemma}.
\end{proof}

\section{Slope independence and the space of photos}
\label{alg}

In \cite{JLM1} and \cite{JLM3}, the second author studied the
\defterm{picture space} of a graph $G$, the algebraic variety of
point-line arrangements in $d$-dimensional space
with an incidence structure given
by $G$.  The rigidity-theoretic behavior of $G$ controls the
geometry of the picture space to a great extent; for instance, the
picture space is irreducible if and only if $G$ is $d$-parallel
independent.

In this section, we study the space $X_{k,d}(M)$ of $(k,d)$-photos,
which is well-defined for any matroid $M$ equipped with a representation.
The photo space plays a role analogous to that of the
picture space of a graph,\footnote{
  The reader should be warned not to take this analogy too literally:
  the picture space of a graph is \emph{not} an instance of the
  photo space of a matroid!}
and the techniques we use to study it are similar to those of \cite{JLM1}. 
In particular, $X_{k,d}(M)$ provides a geometric interpretation of
$m$-Laman independence for all rational numbers $m>1$: it is
equivalent to the space of $(k,d)$-photos.

Throughout this section, we work with a matroid $M$ represented over
a field $\fld$ by nonzero\footnote{
  Our assumption that $M$ contains no loops is purely
  for convenience; trivial (but slightly annoying) modifications
  are necessary when loops are present.}
vectors $v_1,\ldots,v_n \in \fld^r$.  In
addition, let $0<k<d$ be integers, and write $m=\frac{d}{d-k}$.
Recall (Definition~\ref{photo-space}) that the \defterm{space
of $(k,d)$-photos} of $M$ is
  $$\left\{ (\PhMap,W) \in \Hom_\fld(\fld^r,\fld^d) \times \Gr(k,\fld^d)^n:~
    \PhMap(v_i) \in W_i \text{ for all } 1\leq i\leq n \right\}.$$

Note that the photo space is an algebraic subset of
$\Hom_\fld(\fld^r,\fld^d) \times \Gr(k,\fld^d)^n$, hence
a scheme over $\fld$.  The symbol $X_{k,d}(M)$ is a slight
abuse of notation; as defined, the photo space depends on the
representation $\{v_i\}$, and it is not at all clear
to what extent it depends only on the structure of $M$ as an abstract
matroid.  (We will return to this natural question later.)

A key tool in our analysis is a disjoint decomposition of the photo
space into irreducible algebraic subsets called \defterm{cellules}
(in analogy to \cite{JLM1}).  For each photo $(\PhMap,W)$, $\ker \PhMap$
is a linear subspace of $\fld^r$, hence intersects $E$ in some flat $F$ of
$M$.  The idea is to classify photos according to what this flat is.

\begin{defn}
\label{cellule-definition}
For each flat $F \subseteq E$, the corresponding \defterm{cellule} is
  $$
  X_{k,d}^F(M) = \left\{ (\PhMap,W) \in X_{k,d}(M): \ker\PhMap \cap E = F \right\}.
  $$
\end{defn}

By definition, each photo belongs to exactly one cellule; that is,
$X_{k,d}(M)$ decomposes as a disjoint union of the cellules.
Of particular importance are the two extreme cases:

\bigskip\noindent I.
The cellule $X_{k,d}^{\0}(M)$ corresponding to the empty
flat $\0$ is called the \defterm{non-annihilating cellule}.
It is a Zariski open subset of $X_{k,d}(M)$, defined by the conjunction of
open conditions
  \begin{equation} \label{non-ann-conds}
  \PhMap(v_i) \neq 0, \qquad \forall i=1,\dots,n.
  \end{equation}

\bigskip\noindent II.
The cellule $X_{k,d}^E(M)$ corresponding to the improper flat $E$
is called the \defterm{degenerate cellule}.  It is precisely $\{0\} \x
\Gr(k,\fld^d)^n$, where $0$ is the zero map $\fld^r \rightarrow \fld^d$.

\begin{prop} \label{structure-of-cellule}
Let $M$ and $X_{k,d}(M)$ be as above.

\begin{enumerate}

\item[(i)] The natural projection map
    $$X^{\0}_{k,d}(M) \to \Hom_\fld(\fld^r,\fld^d)$$
  gives $X^{\0}_{k,d}(M)$ the structure of an algebraic fiber bundle,
  with fiber $\Gr(k-1,\fld^{d-1})$ and base the Zariski open
  subset of $\Hom_\fld(\fld^r,\fld^d)$ defined by \eqref{non-ann-conds}.  In
  particular, $\dim X^{\0}_{k,d}(M) = dr+n(k-1)(d-k)$.

\item[(ii)] For each flat $F$, $X^F_{k,d}(M) \isom
  X^{\0}_{k,d}(M/F) \times \Gr(k,\fld^d)^F$.  Consequently,
  $X^F_{k,d}(M)$ is an irreducible subvariety of $X_{k,d}(M)$, with
  dimension given by the formula
    \begin{equation} \label{cellule-dim}
    \dim X^F_{k,d}(M) = d(r-r(F)) + (n-|F|)(k-1)(d-k) + |F|k(d-k).
    \end{equation}

\end{enumerate}
\end{prop}

The preceding assertions are more or less immediate from the definition of
cellules and the standard fact that the Grassmannian $\Gr(k,\fld^d)$
has dimension $k(d-k)$.

As in \eqref{projection-map}, let
$\pi$ denote the projection map
  $$
  \Hom_\fld(\fld^r,\fld^d) \times \Gr(k,\fld^d)^n
  \ \overset{\pi}{\longrightarrow} \ \Gr(k,\fld^d)^n,
  $$
and define $M$ to be \defterm{$(k,d)$-slope independent} if
$\pi X^{\0}_{k,d}(M)$ is Zariski dense in $\Gr(k,\fld^d)^n$.
We will denote the Zariski closure of a set $Z$ by $\ov{Z}$.

\begin{thm}
\label{omnibus-theorem}
Let $M$ be a matroid with rank function $r$, represented by vectors
$v_1,\ldots,v_n \in \fld^r$.  Fix positive integers $0<k<d$,
and let $m=\frac{d}{d-k}$.

Then the following are equivalent:
\begin{enumerate}

\item[(i)] $M$ is $(k,d)$-slope independent, i.e., $\pi X^{\0}_{k,d}(M)$ is
dense in $\Gr(k,\fld^d)^n$.

\item[(ii)] $M$ is $m$-Laman independent, i.e., $m \cdot r(F) > |F|$ for every nonempty
flat $F$ of $M$.

\item[(iii)] $\dim X_{k,d}^F(M) < \dim X_{k,d}^{\0}(M)$ for every nonempty flat $F$ of $M$.

\item[(iv)] The photo space $X_{k,d}(M)$ is irreducible.

\item[(v)] The photo space $X_{k,d}(M)$ coincides with the
Zariski closure $\ov{X^{\0}_{k,d}(M)}$ of its non-annihilating cellule.

\end{enumerate}
\end{thm}

\begin{proof}
{\sf (ii)~$\Leftrightarrow$~(iii):}
Compare the cellule dimension formula \eqref{cellule-dim} dimension with
the definition of $m$-Laman independence (Definition~\ref{laman-matroid}).

\smallskip\noindent
{\sf (i)~$\Rightarrow$~(ii):}
For a nonempty flat $F$, write $M|_F$ for the restriction of $M$ to $F$.
Consider the commutative diagram
  \begin{equation}
  \begin{CD}
   X^{\0}_{k,d}(M)    @>{}>>  X^{\0}_{k,d}(M|_F)  \\
  {\pi}{}{\Bigg\downarrow}   &                & {\tilde\pi}{}{\Bigg\downarrow}\\
  \Gr(k,\fld^d)^n    @>{}>> \Gr(k,\fld^d)^F \\
  \end{CD}
  \end{equation}
in which the top horizontal morphism restricts the photo map
$\PhMap$ to the linear span $\fld(F)$ of the vectors in $F$, while
forgetting the $k$-planes $\{W_i\}_{i \in E\sm F}$.
Both vertical
arrows are projections as in \eqref{projection-map}; we use the tilde on the
right-hand map to distinguish them in what follows.
Note that when $\PhMap$ is non-annihilating, its restriction to the span of
$F$ will also be non-annihilating.  Moreover, the bottom horizontal morphism is
surjective.

Now assume that condition (i) holds.  Since the image of
$\pi$ is Zariski dense in the target, so is the image of $\tilde\pi$.
Therefore
  \begin{equation} \label{zariski-inequality}
    d \cdot r(F) + |F|(k-1)(d-k) = \dim X^{\0}_{k,d}(M_F) \geq
    \dim \Gr(k,\fld^d)^F = |F|k(d-k),
  \end{equation}
or in other words, $d \cdot r(F) \geq (d-k)|F|$.
However, scaling a non-annihilating map $\PhMap$
by an element of $\fld^\times$ does not change the line
spanned by any $\PhMap(v_i)$.  Therefore every fiber of
$\tilde\pi$ is at least one-dimensional.
Put differently, when
restricted to $X^{\0}_{k,d}(M|_F)$, the morphism $\tilde\pi$
factors through a $(d \cdot r(F)-1)+|F|(k-1)(d-k)$-dimensional
space of projectivized non-annihilating maps $\PhMap$ in 
$\Pp(\Hom_\fld(\fld(F),\fld^d)$.

Hence, for every nonempty flat $F$, we have the strict inequality
$d\cdot r(F) > (d-k)|F|$, or equivalently
$m \cdot r(F) > |F|$, which is (ii).

\smallskip\noindent
{\sf (iv)~$\Leftrightarrow$~(v):}
Since $X^{\0}_{k,d}(M)$ is Zariski open in $X_{k,d}(M)$,
its closure $\ov{X^{\0}_{k,d}(M)}$ is one of
the irreducible components of $X_{k,d}(M)$.  Thus the full photo
space is irreducible if and only if the non-annihilating
photos are dense.

\smallskip\noindent
{\sf (v)~$\Rightarrow$~(i):}  Suppose that (v) holds.  Then
(i) follows from the observation that
  $$
  \ov{ \pi\left( X^{\0}_{k,d}(M) \right) }
  \supset
  \pi \left(\, \ov{ X^{\0}_{k,d}(M) } \,\right)
    \supset
  \pi( X^E_{k,d}(M) )
  =
  \left( \Pp_\fld^{d-1} \right)^E,
  $$
(the first inclusion is standard, and the second is implied by (v)).

\smallskip\noindent
{\sf (iii)~$\Rightarrow$~(iv):}
We begin by finding an upper bound for
the codimension of every component of the photo space.  Note that
$X_{k,d}(M) = \bigcap_{i=1}^n Z_i$, where
  $$Z_i = \left\{ (\PhMap,W) \in \Hom(\fld^r,\fld^d) \x \Gr(k,\fld^d)^n :~
  \PhMap(v_i)\in W_i \right\}.$$
Let
  \begin{eqnarray*}
    Z'_i  &=& \left\{ (\PhMap,W) \in Z_i :~ \PhMap(v_i)\neq 0 \right\},\\
    Z''_i &=& \left\{ (\PhMap,W) \in Z_i :~ \PhMap(v_i) =   0 \right\}.
  \end{eqnarray*}
Note that $Z'_i$ has codimension $d-k$ in $\Hom_\fld(\fld^r,\fld^d) \x \Gr(k,\fld^d)^n$.
Additionally, $Z''_i$ is contained in the Zariski closure of $Z'_i$,
because the condition $\PhMap(v_i) \in W_i$ (expressed using the Pl\"ucker coordinates
of $W_i$) is satisfied also when $\PhMap(v_i)=0$.  Therefore, every $Z_i$ has
codimension $d-k$, and every irreducible component of their intersection 
$X_{k,d}(M)$ has codimension at most $n(d-k)$.  On the other hand,
by the cellule dimension formula \eqref{cellule-dim}, $n(d-k)$ is precisely
the codimension of the non-annihilating cellule $X^{\0}_{k,d}(M)$.
Hence every irreducible component of $X_{k,d}(M)$ has dimension at least
as large as that of $X^{\0}_{k,d}(M)$.

Now suppose that (iii) holds, so that $\dim X^F_{k,d}(M) < \dim X^\0_{k,d}(M)$
for every $F\neq\0$.  Since the cellules are
all irreducible and disjointly decompose $X_{k,d}(M)$, the irreducible components of 
$X_{k,d}(M)$ must be exactly the closed cellules $\ov{X^F_{k,d}(M)}$
that are contained in the closure of no other cellule.  On the other hand, by the previous
paragraph, every such cellule must have its dimension at least that of
$\dim X^{\0}_{k,d}(M)$, and by (iii) the only possibility is $F=\0$.
Therefore $\ov{X^{\0}_{k,d}(M)}$ is the unique irreducible component.
\end{proof}

The equivalence of (i) and (ii) in Theorem~\ref{omnibus-theorem}
immediately gives the following equality between the slope and Laman complexes.

\begin{cor}
\label{Laman-slope-relation}
Let $m \in \Qq \cap (1,\infty)_\Rr$.  Write $m$ as $\frac{d}{d-k}$, where
$0<k<d$ are integers.

Then $\SSS^{k,d}(M) = \LL^{m}(M)$.
\end{cor}

\begin{rmk} \label{infinite-case}
The condition $d\geq 2$ is implicit in Corollary~\ref{Laman-slope-relation}.
However, there is a sense in which the result is still valid for $d=1$.
Take $k=1$, so that the result asserts that
  $$
  \SSS^{1,d}(M) = \LL^{\frac{d}{d-1}}(M).
  $$
Now, if one establishes conventions properly,
this equality remains valid as $d$ approaches $1$, so that
$m=\frac{d}{d-1}$ approaches infinity.
That is,
  $$
  \SSS^{1,1}(M) = \LL^{\infty}(M) = 2^E.
  $$
Indeed, the full simplex $2^E$ is logically equal to $\SSS^{1,1}(M)$: there is only one 
possible line through any point in $\fld^1$, so the projection map $\pi$ 
is dense. Meanwhile, it is easy to see that $\LL^{\infty}(M)=2^E$, where
we have defined
  $$
  \LL^{\infty}(M):= \lim_{m\rightarrow \infty} \LL^{m}(M).
  $$
\end{rmk}

\begin{rmk} \label{irrational}
For a given matroid $M$ and \emph{irrational} number $m$, it is not hard to see
that there exists a rational number $\tilde m$, chosen sufficiently close to $m$,
such that $\LL^{\tilde m}(M) = \LL^{m}(M)$.  Therefore, 
Corollary~\ref{Laman-slope-relation} actually gives a geometric
interpretation for \emph{every} instance of Laman independence.
\end{rmk}

\begin{rmk} \label{funky-equalities}
Another surprising consequence of Corollary~\ref{Laman-slope-relation} is that
$(k,d)$-slope-independence is invariant under simultaneously scaling $k$ and $d$.
That is, if $\lambda>0$ is an integer, then the Corollary implies that
  $$\SSS^{k,d}(M) = \SSS^{\lambda k,\lambda d}(M).$$
Moreover, if $d$ is divisible by $k$, then $m=d/(d-k)$ is an integer
and $\SSS^{k,d}(M)=\LL^{m}(M)$ is in fact a matroid by Theorem~\ref{Laman-matroidal}~(i).
The geometry behind these phenomena is far from clear.
\end{rmk}

A natural question is to determine the singularities of the photo space.  While we 
cannot do this in general, we can at least say exactly for which matroids $X_{k,d}(M)$
is smooth.  The result and its proof are akin to \cite[Proposition~15]{JLM3},
and do not depend on the parameters $k$ and $d$.

\begin{prop} \label{smooth-Boolean}
Let $M$ be a loopless matroid equipped with a representation $\{v_1,\dots,v_n\}$
as above.  Then, for all integers $0<k<d$, the photo space $X=X_{k,d}(M)$
is smooth if and only if $M$ is Boolean (that is, every ground set element
is an isthmus).
\end{prop}

The assumption of looplessness is harmless, because if $v_i$ is a loop, then
$X_{k,d}(M) \isom \Gr(k,\fld^d) \x X_{k,d}(M\mdel v)$, so $X_{k,d}(M)$
is smooth if and only if $X_{k,d}(M\mdel v)$ is.

\begin{proof}
First, note that the photo space of a direct sum of matroids is precisely the
product of their photo spaces (this can be seen by writing the matrix
for a picture of the direct sum in block-diagonal form).
In particular, if $M$ is Boolean, then
  $$
  X \isom \prod_{i=1}^n \left\{ (\PhMap_i,W_i) \in \fld^d \x \Gr(k,\fld^d):~
  \PhMap_i(v_i) \in W_i \right\},
  $$
and each factor in the product is a copy of the total space of the tautological $k$-plane
bundle over $\Gr(k,\fld^d)$.  In particular, $X$ is smooth.

Now suppose that $M$ is not Boolean; in particular $n > r$.  
Recall
from Proposition~\ref{structure-of-cellule} that the non-annihilating cellule has dimension
$dr+n(k-1)(d-k)$.  Near each non-annihilating photo $\Omega$, the photo space
looks locally like an affine space of this dimension; in particular,
the tangent space $T_\Omega(X)$ has dimension
  \begin{equation} \label{tanOmega}
  \dim T_\Omega(X) = dr+n(k-1)(d-k).
  \end{equation}

Let $\Phi=(\PhMap,W)$ be a ``very degenerate'' photo; that is, $\PhMap = 0$
and all the $k$-planes $W_i$ coincide.  Each $W_i$ can be moved freely throughout the
$i$th Grassmannian, giving $n\cdot\dim\Gr(k,\fld^d) = nk(d-k)$
independent tangent vectors to $X$ at $\Phi$.  On the other hand,
we can also vary the map $\PhMap$ throughout $\Hom(\fld^r,W_i)$,
giving $kr$ more tangent directions that are linearly independent of those
just mentioned.  Therefore
  \begin{equation} \label{tanPhi}
  \dim T_\Phi(X) \geq nk(d-k)+kr.
  \end{equation}
Comparing \eqref{tanOmega} and \eqref{tanPhi}, and doing a little algebra,
we find that
  $$\dim T_\Phi(X) - \dim T_\Omega(X) \geq (d-k)(n-r) > 0.$$
That is, not all points of $X$ have the same tangent space dimension.
Therefore $X$ cannot be smooth.
\end{proof}

\section{Counting photos}
\label{photo-count}

Although it will not be needed in the sequel, we digress to
prove an enumerative result, possibly of independent interest,
about the photo space: when working over a finite field, the
cardinality $|X_{k,d}(M)|$ is an evaluation of the Tutte polynomial
$T_M(x,y)$.

We refer the reader to \cite{BryOx} for details on the
Tutte polynomial.  In what follows, we write $M\mdel v$ and $M/v$ respectively
for the deletion and contraction of $M$ with respect to an element $v$
of its ground set.  We also dispense with the assumption from the
previous section that $M$ contains no loops.  On the other hand,
we add the assumption that the representing vectors $v_1,\ldots,v_n
\in \fld^r$ actually {\it span} $\fld^r$; in other words, $r(M)=r$.
This represents no loss of generality; it is easy to check that
when $r(M)<r$, there is a natural isomorphism
  $$
  X_{k,d}(M) \isom \Hom_\fld( \fld^{r-r(M)},\fld^d) \times X_{k,d}(M'),
  $$
where $M'$ is represented by the same vectors $v_1,\ldots,v_r$, regarded
as elements of the $r(M)$-dimensional subspace of $\fld^r$ that they span.

The following fact \cite[Corollary 6.2.6]{BryOx} is a standard tool for
converting deletion-contraction recurrences to Tutte polynomial evaluations.
We need the \defterm{dual matroid} $M^\perp$, characterized as follows:
when $M$ is represented by the columns $v_1,\ldots,v_n$ of an $r \times n$ matrix of
rank $r$ as above, the dual $M^\perp$ is represented by the columns
$v^*_1,\ldots,v^*_n$ of an $(n-r) \times n$ matrix of rank $n-r$, with
the property that the row space of $M^\perp$
is the nullspace of $M$, and vice versa.  (In purely combinatorial terms,
the bases of $M^\perp$ are the complements of bases of $M$.)

\begin{prop}
\label{Tutte-recursion-converter}
Let $\Psi(M)$ be an invariant of matroids taking values in a commutative
ring $R$, with the following properties:

\begin{tabular}{ll}
{\bf (T1)} \quad & For all matroids $M_1,M_2$,
  $\Psi(M_1 \oplus M_2) = \Psi(M_1) \Psi(M_2)$.\\
{\bf (T2)} & When the ground set of $M$ consists of a single isthmus,
  $\Psi(M) = c$.\\
{\bf (T3)} & When the ground set of $M$ consists of a single loop,
  $\Psi(M) = d$.\\
{\bf (T4)} & When $v$ is neither a loop nor an isthmus of $M$,
  $\Psi(M) = a \Psi(M\mdel v) + b \Psi(M/vf)$.
\end{tabular}

\noindent Then
  $$
  \Psi(M) = a^{r(M^\perp)} b^{r(M)} T_M\left(\frac{c}{b}, \frac{d}{a}\right).
  $$
\end{prop}

Recall \cite[Proposition 1.3.18]{Stanley} that when $\fld$ is a finite field
with $q$ elements, the cardinality of the Grassmannian $\Gr(k,\fld^d)$ 
is given by the {\it $q$-binomial coefficient}
  $$
  \qbin{d}{k} := \frac{[d]!_q}{[k]!_q [d-k]!_q}
  $$
where
  $$
  [n]!_q := [n]_q [n-1]_q \cdots [2]_q [1]_q
  $$
and
  $$
  [n]_q :=\frac{1-q^n}{1-q} = 1 + q + q^2 + \cdots + q^{n-1}.
  $$

We can now state the main result on counting photos.

\begin{thm} \label{Tutte-thm}
Let $\fld$ be the finite field with $q$ elements.  
Let $M$ be a matroid of rank $r$, represented over $\fld$
by vectors $v_1,\dots,v_n$ spanning $\fld^r$, and let $d \geq 2$.  Then
the number of $(k,d)$-photos of $M$ is
$$
|X_{k,d}(M)| \ = \ \qbin{d-1}{k-1}^{r(M^\perp)} \left( q^k \qbin{d-1}{k} \right)^{r(M)}
    T_M \left( \frac{[d]_q}{[d-k]_q} , 
               \frac{[d]]_q}{[k]_q}  \right)
$$
\end{thm}

\begin{proof}
Abbreviate $X_{k,d}(M)$ by $X(M)$, and define $\Psi(M) := |X(M)|$.
We must show that $\Psi$ satisfies the conditions
of Proposition~\ref{Tutte-recursion-converter} with
  $$
  a = \qbin{d-1}{k-1}, \qquad b = q^k\qbin{d-1}{k}, \qquad
  c = q^k\qbin{d}{k},  \qquad d = \qbin{d}{k}.
  $$
(By an easy calculation, the arguments to the Tutte polynomial
in the statement of the theorem are precisely $c/b$ and $d/a$.)

Condition {\bf (T1)} is straightforward.
For {\bf (T2)}, if the ground set of $M$ consists of a single loop,
then $X(M) \isom \Gr(k,\fld^d)$ has cardinality $\qbinsmall{d}{k}$.

If the ground set of $M$ consists of a single isthmus $v$,
then a $(k,d)$-photo of $M$ is just a pair $(\PhMap,W)$ where
$\PhMap: \fld^1 \rightarrow \fld^d$ and $W$ is a $k$-plane
containing $\PhMap(v)$.  Since the image vector $w:=\PhMap(v)$
completely determines the map $\PhMap$, a photo is equivalent to
a pair $(w,W) \in \fld^d \times \Gr(k,\fld^d)$ satisfying
$w \in W$.  Thus the space $X_{k,d}(M)$ is isomorphic to the
tautological $k$-plane bundle over $\Gr(k,\fld^d)$, and its
cardinality is $q^k \qbinsmall{d}{k}$, establishing
condition {\bf (T3)}.

The verification of {\bf (T4)} is the crux of the proof.  If
$v$ is neither a loop nor an isthmus of $M$, we have
the following commutative diagram:

\begin{equation}
\begin{CD}
\label{Jeremy's-diagram}
{\mathcal E}    & \,\, \hookrightarrow \,\, &  X(M)  \\
{\tilde\pi}{}{\Bigg\downarrow}   &                 & {\pi}{}{\Bigg\downarrow}\\
\ov{\mathcal E}        & \,\, \hookrightarrow \,\, &  X(M-v) \\
\end{CD}
\end{equation}

The map $\pi$ sends a $(k,d)$-photo of $M$
to a photo of $M\mdel v$ by forgetting the $k$-plane corresponding
to the vector $v$.  The map $\tilde\pi$ is the restriction of $\pi$ to the
source and target
  $$
  \begin{array}{lll}
  {\mathcal E} &:= \left\{ (\PhMap,W) \in X(M): \phantom{\mdel v}~ \PhMap(v)=0 \right\}
      &~\isom X(M/e) \times \Gr(k,\fld^d)\\
  \ov{\mathcal E} &:=\left\{ (\PhMap,W) \in X(M\mdel v):~ \PhMap(v)=0 \right\}
      &~\isom X(M/e)
  \end{array}
  $$
and corresponds to the projection of $X(M/e) \times \Gr(k,\fld^d)$ onto its first factor.
Meanwhile, the restriction
  $$
  X(M) \sm {\mathcal E}\ \overset{\pi}{\longrightarrow}\ X(M\mdel v) \sm \ov{\mathcal E}
  $$
makes $X(M) \sm {\mathcal E}$ into a bundle with
fiber $\Gr(k-1,\fld^{d-1})$.  Consequently
  \begin{align*}
  \left| X(M) \sm {\mathcal E} \right|
          &= \qbin{d-1}{k-1} \left| X(M\mdel v) \sm \ov{\mathcal E}\right| \\
\intertext{and}
  \Psi(M) &= \left|{\mathcal E}\right|
           + \qbin{k-1}{d-1} \left( \Psi(M\mdel v) - \left|\ov{\mathcal E}\right| \right) \\
          &= \qbin{d-1}{k-1} \Psi(M\mdel v)
           + \qbin{d}{k} \Psi(M/v) - \qbin{d-1}{k-1}\Psi(M/v) \\
          &= \qbin{d-1}{k-1} \Psi(M\mdel v) + q^k \qbin{d-1}{k} \Psi(M/v)
  \end{align*}
where the last equality uses the $q$-Pascal recurrence 
\cite[Chapter 1, \S1.3, Equation (17b)]{Stanley}
$$
\qbin{d}{k} = q^k \qbin{d-1}{k} + \qbin{d-1}{k-1}.
$$
\end{proof}

Since the Tutte polynomial of $M$ does not depend on the choice
of representation, neither does the number of photos.
Theorem~\ref{Tutte-thm} also implies a curious symmetry
between the number of photos of a matroid $M$ and of its dual $M^\perp$.
Since $T_{M\!^\perp}(x,y) = T_M(y,x)$ \cite[Prop.~6.2.4]{BryOx} and
$\qbinsmall{d}{k}=\qbinsmall{d}{d-k}$, we have:

\begin{cor}
\label{curious-Tutte-symmetry}
Let $M$ and $M^\perp$ be dual represented matroids.  Then
$$
q^{d \cdot r(M)} |X_{d-k,d} (M^\perp)| = q^{(d-k)n} |X_{k,d}(M)|.
$$
\end{cor}

It would be nice to have a more direct explanation for Corollary~\ref{curious-Tutte-symmetry}.

\begin{rmk}
A topological commutative diagram analogous to \eqref{Jeremy's-diagram}
was exploited by the second author in \cite{JLM2} to compute the Poincar\'e series
of picture spaces of graphs over $\Cc$ as an analogous Tutte polynomial evaluation.  
In contrast, when $\fld=\Rr$ or $\Cc$, the topology of the photo space
is much simpler.  Indeed, there is a deformation retraction of
$X_{k,d}(M)$ onto its degenerate cellule:
    $$\begin{array}{cccc}
      F: & [0,1] \x X_{k,d}(M) &\to& X^E_{k,d}(M) \\
         & (\lambda,(\PhMap,W)) &\mapsto& (\lambda\PhMap,W).
    \end{array}$$
Hence $X_{k,d}(M)$ is homotopy equivalent to the degenerate cellule $X^{\0}_{k,d}(M)$,
which is homeomorphic to $\Gr(k,\fld^d)^n$ (see Definition~\ref{cellule-definition}).
\end{rmk}

\section{Rigidity and parallel independence}
\label{rigidity-section}

In this section, we examine more closely the special cases $k=1$ and $k=d-1$ of
$(k,d)$-slope independence for a represented matroid $M$.  It turns out that
they are intimately related to the $d$-dimensional generic rigidity 
matroid $\R^d(M)$ and the $d$-dimensional generic 
hyperplane-marking matroid $\HH^d(M)$.
Throughout the section, let $M$ be a matroid represented by vectors
$E=\{v_1,\ldots,v_n\}$ spanning $\fld^r$, and let $d>0$ be an integer. 

\subsection{Interpreting $\R^d(M)$ and $\HH^d(M)$}

Recall (Definition~\ref{rigidity-matroid}) that the
\defterm{$d$-dimensional rigidity matroid}
is represented over $\fld(\PhMap)$ by the vectors
  $$
  \{v_i \otimes \PhMap(v_i)\}_{i=1}^n \subset \fld^r \otimes_\fld \fld(\PhMap)^d.
  $$
where $\fld(\PhMap)$ is the extension of $\fld$ by $dr$ transcendentals
(the entries of the matrix $\PhMap: \fld^r \rightarrow \fld(\PhMap)^d$).
The complex $\R^d(M)$ is defined to be the complex of independent sets of this matroid.
The \defterm{$d$-rigidity matrix} $R^d(M)$ is the $n \times dr$ matrix whose rows represent
$\R^d(M)$.

Recall also (Definition~\ref{hyperplane-marking-matroid}) that the
\defterm{$d$-dimensional hyperplane-marking matroid}
is represented over $\fld(\PhMap,n)$ by the vectors
  $$
  \{v_i \otimes \normal_i)\}_{i=1}^n \subset \fld^r \otimes_\fld \fld(\PhMap,\normal)^d.
  $$
where $\fld(\PhMap)$ is the extension of $\fld$ by $dr+(d-1)n$ transcendentals
(the $dr$ entries of the matrix $\PhMap$, and the $(d-1)n$ coordinates of
the normal vectors $\normal_i$ to $\PhMap(v_i)$).
The complex $\HH^d(M)$ is defined to be the complex of independent sets of this matroid.
Denote by $H^d(M)$ the $n \times dr$ matrix whose rows represent $\HH^d(M)$.

To interpret $R^d(M)$ and $H^d(M)$, we study their (right) nullspaces.
Both matrices have row vectors in $\fld^r \otimes_\fld \fld^d$,
so their nullvectors live in the same space.  It will be convenient to
freely use the identifications
  $$
  \fld^r \otimes_\fld \fld^d \isom
  (\fld^r)^* \otimes_\fld \fld^d \isom 
  \Hom_\fld(\fld^r,\fld^d).
  $$
The second of these isomorphisms is canonical; the first comes from
identifying $\fld^r$ and $(\fld^r)^*$ by the standard bilinear form on $\fld^r$,
  $$\langle x,y \rangle =\sum_{i=1}^r x_i y_i,$$
whose associated quadratic form is
  $$Q(x) = \langle x,x \rangle = \sum_{i=1}^r x_i^2.$$

With these identifications, for every $\psi \in \fld^r \otimes_\fld \fld^d \isom
\Hom_\fld(\fld^r,\fld^d)$, $v \in \fld^r$, and $x \in \fld^d$, the induced bilinear form
on $\fld^r \otimes_\fld \fld^d$ has the property
  $$
  \langle v \otimes x, \psi \rangle = \langle x, \psi(v) \rangle.
  $$

\begin{prop}
\label{rigidity-parallel-interpretations}
Let $M$ be a matroid represented by $E$ as above, and
let $\psi \in \fld^r \otimes_\fld \fld^d \isom \Hom_\fld(\fld^r,\fld^d)$.

\begin{enumerate}

\item[(i)] The vector $\psi$ lies in $\ker H^d(M)$ if and only if
$(\PhMap+\psi)(v_i)$ is normal to $\normal_i$
for all $i=1,2,\ldots,n$.  

{\rm (In other words, the nullspace of $H^d(M)$ is the space of directions in
which one can modify the map $\PhMap$ while keeping the image of $v_i$ lying on the
same hyperplane normal to $\normal_i$ for each $i$.)}

\item[(ii)] Provided that $\fld$ does not have characteristic $2$, the vector $\psi$ lies in
$\ker R^d(M)$ if and only if
  $$
  Q\big( (\PhMap + \epsilon\psi)(v_i) \big) \equiv Q\big( \PhMap(v_i) \big) \mod \epsilon^2
  $$
for each $i=1,2,\ldots,n$.  

{\rm (In other words, the nullspace of $R^d(M)$ is the space of
infinitesimal modifications one can make to $\PhMap$ while keeping the values
of the quadratic form $Q$ on the images of the $v_i$ constant (up to first order)
for each $i$.)}
\end{enumerate}
\end{prop}

\begin{proof}
For (i), note that
\begin{align*}
\langle \normal_i, (\PhMap+\psi)(v_i) \rangle = 0 \
&\iff\ \langle \normal_i, \PhMap(v_i) \rangle + \langle \normal_i, \psi(v_i) \rangle = 0\\
&\iff\ \langle \normal_i, \psi(v_i) \rangle = 0\\
&\iff\ \langle v_i \otimes \normal_i, \psi  \rangle = 0.
\end{align*}

For (ii), the expression
$$
Q( (\PhMap + \epsilon \psi)(v_i) ) =
Q( (\PhMap(v_i) ) + 2 \epsilon \langle \PhMap(v_i), \psi(v_i) \rangle + \epsilon^2 Q( \psi(v_i) )
$$
is congruent to $Q( \PhMap(v_i) )$ modulo $\epsilon^2$ if and only
if $\langle \PhMap(v_i), \psi(v_i) \rangle =0$ (since $\fld$ does not have characteristic $2$).
But $\langle \PhMap(v_i), \psi(v_i) \rangle = \langle v_i \otimes \PhMap(v_i), \psi \rangle$,
completing the proof.
\end{proof}

\begin{rmk}
\label{Jacobian-interpretation}
Part (i) of Proposition~\ref{rigidity-parallel-interpretations}
is a rephrasing of the following familiar fact from rigidity theory:
the rigidity matrix $R^d(M)$ may be regarded as the Jacobian matrix
(after scaling by $\frac{1}{2}$) of the map 
$$
\begin{array}{lll}
\Hom_\fld(\fld^r,\fld^d) & \longrightarrow & \fld^n \\
\PhMap                   & \longmapsto     & ( Q(\PhMap(v_i) )_{i=1}^n\ .\\
\end{array}
$$
\end{rmk}

The following instance of the hyperplane-marking matroid generalizes
the notion of the $d$-parallel matroid of a graph (see
\eqref{graph-parallel-vectors}).
Denote by $(d-1)M$ the matroid whose ground set consists of
$d-1$ copies of each vector in $E$.
The \defterm{$d$-parallel matrix} of $M$ is defined as
$H^d((d-1)M)$, and the matroid represented by its rows
is the \defterm{($d$-dimensional, generic) parallel matroid}
$\PP^d(M) := \HH^d((d-1)M)$.
Part (ii) of Proposition~\ref{rigidity-parallel-interpretations}
leads to an interpretation of the geometric meaning carried by the
$d$-parallel matrix:

\begin{cor}
Let $\psi \in \fld^r \otimes_\fld \fld^d \isom \Hom_\fld(\fld^r,\fld^d)$.
Then $\psi \in \ker P^d(M)$ if and only if $(\PhMap+\psi)(v_i)$ is parallel
to $\PhMap(v_i)$ for all $i=1,2,\ldots,n$.  
\end{cor}

\begin{proof}
Since there are $d-1$ copies of the vector $v_i$ in $(d-1)M$, there will be
$(d-1)$ accompanying normal vectors to $\PhMap(v_i)$.  Because these normals are
chosen with generic coordinates, the only vectors normal to all $d-1$ of them
are those parallel to $\PhMap(v_i)$.  Now apply
Proposition~\ref{rigidity-parallel-interpretations}.
\end{proof}

As in classical rigidity theory, both $R^d(M)$ and $H^d(M)$ have certain obvious nullvectors.

\begin{prop}
\label{obvious-nullvectors}
Let $\psi \in \fld^r \otimes_\fld \fld^d \isom \Hom_\fld(\fld^r,\fld^d)$.
\begin{enumerate}
\item[(i)] Given any skew-symmetric $d\times d$ matrix $\sigma \in \fld^{d \times d}$,
the map $\sigma \circ \psi$, when identified with a vector in $\fld^r \otimes \fld^d$,
lies in the nullspace of $R^d(M)$.
\item[(ii)] The map $\psi$, when identified with a vector in $\fld^r \otimes \fld^d$,
lies in the nullspace of $H^d(M)$.
\end{enumerate}
\end{prop}

\begin{proof}
Assertion (ii) is immediate from the interpretation of the nullspace
of $H^d(M)$ given in Proposition~\ref{rigidity-parallel-interpretations}.

To prove (i), we define
  $$S := \Zz[\PhMap, \sigma, v]/\left(\sigma_{ji}=-\sigma_{ij}\right),$$
the polynomial ring in the entries of the matrices
$\PhMap, \sigma, v_1,\ldots,v_n$.  We wish to show that
  \begin{equation}
  \label{desired-Jacobian-identity}
  R^d(M) (\sigma \circ \PhMap) = 0
  \end{equation}
in $S$.  In fact, we will show by a formal calculation that
$2R^d(M) (\sigma \circ \PhMap) = 0$.  Since $2$ is a non-zero-divisor
in $S$, this will imply that \eqref{desired-Jacobian-identity} holds in
$S$, hence remains valid when we pass to $S \otimes_{\Zz} \fld$ and
specialize the entries of $v_1,\ldots,v_n,\sigma$ to elements of $\fld$.

The calculation\footnote{
  This calculation is identical to that usually used to show that
  the orthogonal group with respect to the quadratic form $Q$ on $\fld^d$
  has its Lie algebra equal to the space of skew-symmetric matrices.}
actually takes place in $S[\epsilon]/(\epsilon^2)$.
Since $\sigma^T=-\sigma$, one has for all $x\in\fld^d$
\begin{align*}
Q((I_d+ \epsilon \sigma)(x)) &= Q(x)+ \epsilon \langle x, \sigma(x) \rangle +
                                       \epsilon \langle \sigma(x), x \rangle + 
                                       \epsilon^2 Q( \sigma(x)) \\
                             &= Q(x)+ \epsilon \left( \langle x, \sigma(x) \rangle +
                                                \langle x, \sigma^T(x) \rangle \right) + 
                                       \epsilon^2 Q( \sigma(x)) \\
                             &\equiv Q(x) \mod \epsilon^2
\end{align*}

Taking $x=\PhMap(v_i)$, the function $f$ defined by
$f(\PhMap) := Q(\PhMap(v_i))$ has the property
  $$
  f(\PhMap+ \epsilon \sigma \circ \PhMap) \equiv f(\PhMap) \mod \epsilon^2.
  $$
On the other hand, expanding $f$ as a Taylor polynomial yields
  $$
  f(\PhMap+ \epsilon \sigma \circ \PhMap) \equiv
    f(\PhMap) + \epsilon \langle \gradient_\PhMap(f), \sigma \circ \PhMap \rangle 
         \mod \epsilon^2.
  $$
where $\gradient_\PhMap(f)$ is the gradient of $f$ with respect to
the entries of $\PhMap$.
Therefore $\langle \gradient_\PhMap(f), \sigma \circ \PhMap \rangle = 0$.
On the other hand, by Remark~\ref{Jacobian-interpretation}, 
the $i^{th}$ row of $R^d(M)$
is exactly $\half\gradient_\PhMap(f)$.  So $2R^d(M) \sigma \circ \PhMap = 0$
as desired.
\end{proof}

\subsection{The Nesting Theorem}

We have arrived at one of the main results of the paper,
the \emph{Nesting Theorem},
which explains the relationship between the various independence
systems associated to an arbitrary representable matroid $M$.
In the special case that $M$ is graphic and the ambient dimension $d$ is
2, the Nesting Theorem gives what we have called the {\it planar trinity}
(Corollary~\ref{planar-case} below).

\begin{thm}[The Nesting Theorem]
\label{Laman-rigidity-relation}
Let $M$ be a matroid represented by vectors
$E=\{v_1,\ldots,v_n\} \subset \fld^r$, and let $d>1$ be an integer. 
Then
  $$
  \SSS^{1,d}(M) \subseteq \DRig(M) \subseteq \DLam(M) = \HH^d(M) \quad (= \SSS^{d-1,d}(M)).
  $$
\end{thm}

\begin{proof}
We first prove that $\DRig(M) \subseteq \DLam(M)$.
It suffices to show that whenever $d \cdot r(M) \leq n$, there is
an $\fld(\PhMap)$-linear dependence among the vectors
  $$\{v_i \otimes \PhMap(v_i)\}_{i=1}^n \subset \fld^r \otimes_\fld \fld^d$$
that form the $n$ rows of $R^d(M)$.
Since $E$ spans a subspace of $\fld^r$ isomorphic
to $\fld^{r(M)}$, the rows of $R^d(M)$ actually lie in a subspace of
dimension $d \cdot r(M)$.  If $d \cdot r(M) < n$, then the
desired linear dependence is immediate.  On the other hand, if
$d \cdot r(M) = n$, then Proposition~\ref{obvious-nullvectors}
implies that the rows of $R^d(M)$ lie in a \emph{proper} subspace of 
$\fld^{r(M)} \otimes \fld^d$, hence are linearly dependent.

If we replace $v_i \otimes \PhMap(v_i)$ with $v_i \otimes \normal_i$,
the same argument shows that $\HH^d(M) \subseteq \DLam(M)$.

\smallskip

Next we prove that $\SSS^{1,d}(M) \subseteq \DRig(M)$.
Assume that the rows of $R^d(M)$
are dependent; we will show that $M$ is $(k,d)$-slope
dependent for $k=1$.

We begin with the observation that 
  $$
  \SSS^{k,d}(M) = \LL^{\frac{d}{d-k}}(M) \subseteq \LL^d(M).
  $$
The equality is Corollary~\ref{Laman-slope-relation}, and the
inclusion follows from the definition of $\LL^m(M)$ (because $\frac{d}{d-k} \leq d$).
In particular, if $M$ is $d$-Laman dependent then $M$ is
automatically $(k,d)$-slope dependent; we may therefore assume 
that $M$ is $d$-Laman independent.  Without loss of
generality, $d \cdot r(M) \geq n$, so the dependence of
the rows of $R^d(M)$
implies the vanishing of \emph{every} one of its $n \times n$ minor subdeterminants.
Moreover, by Theorem~\ref{Edmonds-Laman-equiv},
$M$ admits a $d$-Edmonds decomposition (see Definition~\ref{Edmonds-decomp}).
Associating the vectors $v_1,\dots,v_n$
with their indices $[n]=\{1,\dots,n\}$, we may write this Edmonds decomposition
concisely as $[n] = \bigsqcup_{j=1}^d I_j$.

\smallskip
\begin{quote}
{\bf Claim.}
There exists an $n \times n$ minor $\xi$ of $R^d(M)$
that is a nonzero multihomogeneous polynomial in the coordinates of the
vectors $\PhMap(v_i)$.
\end{quote}
\smallskip

\noindent Given the claim, if $\xi$
vanishes on the non-annihilating cellule $X^{\0}_{k,d}(M)$
of the photo space, then the projection on
$X^{\0}_{k,d}(M) \rightarrow \Gr(k,\fld^d)$ is not Zariski dense,
because the homogeneous coordinates of the $\PhMap(v_i)$
are in fact the Pl\"ucker coordinates on $\Gr(k,\fld^d)$.
Hence by Theorem~\ref{omnibus-theorem}, the claim is
all we need for the present theorem.

Let $x^{(i)}:=\PhMap(v_i)$, and let $v_i=[ v_{i1} \cdots v_{ir} ]^T$.
Group the columns of $R=R^d(M)$
in blocks, so that the $i^{th}$ row of $R$ is
  $$
  \big[ v_{i1} x^{(i)}_1 \, \cdots\, v_{ir} x^{(i)}_1 \ \big\vert \ 
        v_{i1} x^{(i)}_2 \, \cdots\, v_{ir} x^{(i)}_2 \ \big\vert \
        \cdots \ \big\vert \
        v_{i1} x^{(i)}_d \, \cdots\, v_{ir} x^{(i)}_d \big].
  $$
Each $n \times n$ submatrix $R_A$ of $R$ is indexed by some choice of an $n$-element
subset $A$ of the $dr$ columns.  Letting $A_i$ be the subset of $A$ coming from 
columns in the $i^{th}$ block, one obtains a sequence of subsets 
$A_1,\ldots,A_d \subset [r]$ with $n=|A| = \sum_{j=1}^d |A_j|$.
Then
  $$
  \det R_A ~=~
    \sum_I \sgn(I) \sum_{\sigma_1,\dots,\sigma_d} \sgn(\sigma_1) \dots \sgn(\sigma_d)
      \prod_{j=1}^d \prod_{i \in I_j} v_{i,\sigma_j(i)} x^{(i)}_j.
  $$
Here the first sum ranges over all partitions $I = \{I_1,\dots,I_d\}$ of $[n]$ with $d$
parts, the second sum ranges over all $d$-tuples of bijections $\sigma_j:\ I_j\to A_j$,
and $\sgn(C),\sgn(\sigma_j) \in \{\pm 1\}$ (there are explicit formulas for these signs,
but we won't need them).  This expression may be simplified:
  \begin{eqnarray*}
   \det R_A 
    &=& \sum_I \sgn(I) \left( \prod_{j=1}^d \prod_{i \in I_j} x^{(i)}_j \right)
         \left( \prod_{j=1}^d \ \  \sum_{\sigma_j: I_j \to A_j}
             \sgn(\sigma_j) v_{i,\sigma_j(i)} \right) \\
    &=& \sum_I \left( \prod_{j=1}^d \prod_{i \in I_j} x^{(i)}_j \right)
         \left( \sgn(I) \prod_{j=1}^d  \det V_{I_j,A_j} \right)
  \end{eqnarray*}
where $V_{I_j,A_j}$ is the submatrix of 
$[v_{ik}]_{i=1,\dots,n, k=1,\dots,r}$ with rows $I_j$ and columns $A_j$.
Note that $\det( V_{I_j,A_j} ) \in \fld$, so the calculation implies that
$\det R_A$ is a multihomogeneous polynomial
in the coordinates $\{ x^{(i)}_j \}$ with coefficients in $\fld$.

By the definition of an Edmonds decomposition,
the sets $I_1,\dots,I_d$ are independent in $M$.
Hence there is some subset $A_j \subseteq [r]$
with $\det V_{I_j,A_j} \neq 0$.  The monomial
corresponding to this choice of $I_j$'s and $A_j$'s has a nonzero coefficient
in the multihomogeneous polynomial $\xi=\det R_A$.  Therefore $\xi\neq 0$,
establishing the claim and completing the proof that
$\SSS^{1,d}(M) \subseteq \DRig(M)$.

Replacing $R^d(M)$ with $H^d(M)$, $k=1$ with $k=d-1$,
and $\PhMap(v_i)$ with $\normal_i$ throughout, the same
argument shows that $\SSS^{d-1,d} \subseteq \HH^d(M)$.
Since $\SSS^{d-1,d}(M) = \DLam(M)$ by
Corollary~\ref{Laman-slope-relation},
we are done.
\end{proof}

The case $d=2$ is very special.  Recall that $\PP^d(M)=\HH^d((d-1)M)$,
so $\PP^2(M)=\HH^2(M)$.  Indeed, the Nesting Theorem implies much
more:

\begin{cor}
\label{planar-case}
Let $M$ be a matroid represented as above.  Then
  $$
  \SSS^{1,2}(M) = \R^2(M)  = \LL^2(M) = \HH^2(M) = \PP^2(M).
  $$
\end{cor}


\begin{rmk}
Setting $d=1$ collapses the Nesting Theorem to
  $$
  \SSS^{k,\infty}(M) =  \R^1(M) = \LL^1(M) = M.
  $$
However, these phenomena are somewhat more trivial.  To make sense of the
complexes $\SSS^{k,\infty}(M)$ and $\LL^1(M)$, consider the identity
$\SSS^{k,d}(M) = \LL^{\frac{d}{d-k}}(M)$ of Corollary~\ref{Laman-slope-relation}.
Fixing $k$ and letting $d\to\infty$ (as a positive integer), we obtain
$\SSS^{k,\infty}(M) =  \R^1(M) = \LL^1(M)$.  On the other hand,
it is an easy consequence of the definitions of $\LL^m(N)$ and $\R^d(M)$ that
$\lim_{m \rightarrow 1^+} \LL^m(M)=M=\R^1(M)$.
\end{rmk}

\begin{rmk}
There is in fact a simple explicit isomorphism between
the matroids $\R^2(M)$ and $\HH^2(M)$ ($=\PP^2(M)$).
Let $\rho$ be  the ``$\pi/2$ rotation'' $\fld^2 \rightarrow \fld^2$ given by 
  $$
  \left[ \begin{matrix} 
  0 & -1 \\
  1 & 0 
  \end{matrix} \right].
  $$
Then $\rho( \PhMap(v_i) ) = \normal_i$, a generic normal to the
generic image vector $\PhMap(v_i)$, and the invertible linear operator
$1_{\fld^r} \otimes \rho$ on $\fld^r \otimes_\fld \fld^d$
sends $v_i \otimes \PhMap(v_i)$ to $v_i \otimes \normal_i$.
\end{rmk}

\begin{rmk}
When $d \geq 3$, the inclusion $\R^d(M) \subseteq \LL^d(M)$
is usually strict.
By Proposition~\ref{obvious-nullvectors}, the nullspace of
$R^d(M)$ contains the $\binom{d}{2}$-dimensional space of all
vectors of the form $\sigma \circ \PhMap$, as $\sigma$ ranges over all
skew-symmetric matrices in $\fld^{d \times d}$.  Consequently,
every $d$-rigidity-independent subset $A \subseteq E$ must satisfy
$|A| \leq d \cdot r(A) - \binom{d}{2}$.  On the other hand,
there may exist $d$-Laman independent sets $A$ of cardinality up to
$d \cdot r(A)-1$.
\end{rmk}

\section{Examples: Uniform matroids}
\label{examples-section}

Let $E$ be a ground set with $n$ elements.  The \emph{uniform
matroid} of rank $r$ on $E$ is defined to be the
matroid whose independent sets are
  $$\Umat{r}{n} = \{F \subseteq E:~ |F|\leq r\}.$$
Broadly speaking, $\Umat{r}{n}$ can be regarded as the
matroid represented by $n$ generically chosen vectors in
$\fld^r$, where $\fld$ is a sufficiently large field.

Predictably, the $d$-Laman independence complex on $\Umat{r}{n}$
is also a uniform matroid for every $d$.  More surprising is that
$d$-Laman independence carries nontrivial geometric information
about sets of $n$ generic vectors in $r$-space---specifically
coplanarity for $\Umat{2}{3}$ and the cross-ratio for $\Umat{2}{4}$.

\begin{prop} \label{uniform}
Let $\Umat{r}{n}$ be the uniform matroid of rank $r$
on $n$ elements, and let $d\in(1,\infty)_\Rr$.
Then
  \begin{equation} \label{uniform-Laman}
  \LL^d(\Umat{r}{n}) = \Umat{s}{n} \qquad \text{where } s = \min(\lceil dr-1 \rceil, n).
  \end{equation}
and
  \begin{equation} \label{uniform-slope}
  \SSS^{k,d}(\Umat{r}{n}) = \Umat{t}{n}
    \qquad \text{where } t = \min\left( \lceil\tfrac{dr}{d-k}-1\rceil, n \right).
  \end{equation}
\end{prop}

\begin{proof} We know that $\LL^d(\Umat{r}{n})$ is a simplicial complex, and
it is easy to see that the criteria for $F$ to be $d$-Laman independent
can depend only depend on the cardinality $|F|$.  Therefore
  \begin{eqnarray*}
  \LL^d(\Umat{r}{n})
  &=& \{F \subseteq E:~ d\cdot r(F') > |F'| \quad \text{for all nonempty } F'\subseteq F\} \\
  &=& \{F \subseteq E:~ d\cdot r(F) > |F|\} \\
  &=& \{F \subseteq E:~ |F| < dr\}\\
  &=& \Umat{s}{n},
  \end{eqnarray*}
which is \eqref{uniform-Laman}.  Applying Corollary~\ref{Laman-slope-relation}
to \eqref{uniform-Laman} gives \eqref{uniform-slope}.
\end{proof}

\begin{ex}[$\Umat{2}{3}$]
\label{U23-example}
Let $\fld$ be any field, and let $e_1, e_2$ 
be the standard basis vectors in $\fld^2$.
The matroid $M = \Umat{2}{3}$ is represented by the vectors
$\{e_1, e_1+e_2, e_2\} \subset \fld^2$; this representation is
unique up to the action of the projective general linear group.
By Proposition~\ref{uniform},
  $$
  \LL^d(\Umat{2}{3}) = 
  \begin{cases}
  \Umat{2}{3}       & \text{ if } d \in (1, \frac{3}{2}]_\Rr \\
  \Umat{3}{3}       & \text{ if } d \in (\frac{3}{2}, \infty)_\Rr
  \end{cases}
  \qquad\text{and}\qquad
  \SSS^{1,d}(\Umat{2}{3}) = 
  \begin{cases}
  \Umat{3}{3}       & \text{ if } d = 2 \\
  \Umat{2}{3}       & \text{ if } d \in \{3,4,\ldots \}.
  \end{cases}
  $$
We now consider what these equalities mean in terms of slopes.
Let $\PhMap: \fld^2 \to \fld^d$ be a linear transformation.
If $d=2$, then the images $\PhMap(e_1), \PhMap(e_1+e_2), \PhMap(e_2)$
can have arbitrary slopes as $\PhMap$ varies.  This is why
$\SSS^{1,2}(\Umat{2}{3}) = \Umat{3}{3}$.  On the other hand,
when $d \geq 3$, those three vectors must be coplanar.  This imposes
a nontrivial constraint on the homogeneous coordinates for the
lines spanned by the three images, and explains why
$\SSS^{1,d}(\Umat{2}{3}) = \Umat{2}{3}$.

By direct calculation, the vectors
  $$
  e_1 \otimes \PhMap(e_1), \quad (e_1+e_2) \otimes \PhMap(e_1+e_2), \quad
  e_2 \otimes \PhMap(e_2)
  $$ 
are linearly dependent if and only if $d=1$.  Therefore
  $$
  \DRig(\Umat{2}{3}) = 
  \begin{cases}
  \Umat{2}{3}       & \text{ if } d = 1, \\
  \Umat{3}{3}       & \text{ if } d \in \{2,3,\ldots\}.
  \end{cases}
  $$
In this case, the inclusions $\R^d(M) \subseteq \LL^d(M)$ 
given by Theorem~\ref{Laman-rigidity-relation} turn out to be equalities.
\end{ex}

\begin{ex}[$\Umat{2}{4}$]
\label{U24-example}
Let $\fld$ be a field of cardinality $>2$,
let $\mu\in\fld\sm\{0,1\}$,
and let $e_1,e_2$ be the standard basis vectors in $\fld^2$.
The four vectors
  $$
  \{ e_1, \ \ e_1+e_2, \ \ e_2, \ \ e_1+\mu e_2 \},
  $$
represent $M=\Umat{2}{4}$ over $\fld$.  Again,
this representation is unique up to projective equivalence.
By Proposition~\ref{uniform},
  $$
  \LL^d(\Umat{2}{4}) = 
  \begin{cases}
  \Umat{2}{4} & \text{ if } d \in (1, \frac{3}{2}]_\Rr \\
  \Umat{3}{4} & \text{ if } d \in (\frac{3}{2},2]_\Rr \\
  \Umat{4}{4} & \text{ if } d \in (2, \infty)_\Rr
  \end{cases}
  \qquad\text{and}\qquad
  \SSS^{1,d}(\Umat{2}{4}) = 
  \begin{cases}
  \Umat{3}{4} & \text{ if } d = 2 \\
  \Umat{2}{4} & \text{ if } d \in \{3,4,\ldots \}.
  \end{cases}
  $$

Why is this correct from the point of view of slopes?
From Example~\ref{U23-example}, we know that when $d \geq 3$, the lines spanned by
the images of any three of the four vectors must be coplanar, so there is an algebraic
dependence among the homogeneous coordinates for these three lines.
For $d=2$, this does not happen;
the slopes of the images of any triple can be made arbitrary.  However,
applying a linear transformation to the representing vectors does not
change their \emph{cross-ratio} (in this case $\mu$), so the fourth image vector
is determined by the first three.  This is the geometric interpretation
of the combinatorial identity $\SSS^{1,2}(\Umat{2}{4}) = \Umat{3}{4}$.

Direct calculation shows that every three of the four vectors
  $$
  w_1 := e_1 \otimes \PhMap(e_1), \quad w_2 := (e_1+e_2) \otimes \PhMap(e_1+e_2), \quad
  w_3 := e_2 \otimes \PhMap(e_2), \quad w_4 := (e_1+ \mu e_2) \otimes \PhMap(e_1 + \mu e_2)
  $$ 
are linearly dependent when $d=1$, but independent for all $d\geq 2$.
When $d\geq 2$, there is an additional, less obvious linear dependence:
$(\mu-1) w_1 - \mu w_2 + (\mu-\mu^2) w_3 + w_4 = 0$.  Consequently
  $$
  \DRig(\Umat{2}{4}) = 
  \begin{cases}
  \Umat{2}{4}       & \text{ if } d = 1, \\
  \Umat{3}{4}       & \text{ if } d \in \{2,3,\ldots \}.
  \end{cases}
  $$
This calculation is independent of the particular 
coordinates chosen for the representing vectors, even up to
projective equivalence (that is, up to the choice of the parameter $\mu$):
that is, $\DRig(\Umat{2}{4})$ is a \emph{combinatorial} invariant.

On the other hand, unlike the situation for $\Umat{2}{3}$,
the inclusions $\R^d(M) \subseteq \LL^d(M)$ given by
Theorem~\ref{Laman-rigidity-relation} turn out to be strict.
In particular, $\R^\infty(M)$ is not Boolean while $\LL^\infty(M)$ is always
Boolean.  This behavior deviates notably from the case of graphic 
matroids (see Proposition~\ref{graphic-boolean-prop} below).
\end{ex}

\section{More on $\DRig(M)$: invariance and stabilization}
\label{invariance-section}

The examples in the previous section raise some natural questions.
Clearly $\LL^m(M)$ is a {\it combinatorial invariant} of $M$,
that is, it does not depend on the choice of representation, nor the field
of representation.  Hence by Corollary~\ref{Laman-slope-relation},
the same is true for $\SSS^{k,d}(M)$, and in particular $\HH^d(M)$ and $\PP^d(M)$.
But what about $\DRig(M)$?  This is an issue which does not arise in
classical rigidity theory, as the graphic matroid $M(G)$ is always represented
by the vectors
  \begin{equation} \label{std-graphic}
  \{e_i-e_j: \{i,j\} \in E(G) \}
  \end{equation}
where $e_i$ is the $i^{th}$ standard basis vector in $\Rr^{|V(G)|}$.
In fact, Proposition~\ref{projective-invariance} below will show that
$\DRig(M)$ is a projective invariant of a matroid represented over a given field.
A result of N. White shows that graphic matroids, and more generally matroids
that can be represented over $\fld_2$, are projectively unique when represented over any fixed field;
see, e.g., \cite[Proposition 1.2.5]{White}.
It will follow that $\DRig(M(G))$ is a combinatorial invariant of a graphic matroid $M(G)$
over any fixed field. 

We begin by recalling the notion of projective equivalence for representations
of a matroid. Two sets of vectors $E=\{v_1,\ldots,v_n\}$, $E'=\{v'_1,\ldots,v'_n\} \subset \fld^r$
are called \defterm{projectively equivalent} if there are nonzero scalars
$c_1,\dots,c_n \in \fld^\x$ and an invertible linear
transformation $g \in GL_r(\fld)$, such that $v'_i=g(c_iv_i)$ for every $i$.
It is easy to see that in this case, the matroids represented by $E$ and $E'$
are combinatorially identical.  As we now show, the same is true for their
$d$-rigidity matroids.

\begin{prop}\label{projective-invariance}
Let $M,M'$ be matroids represented by projectively equivalent sets
$E,E' \subset \fld^r$, and let $d\geq 2$.  Then $\DRig(M)=\DRig(M')$.
\end{prop}

\begin{proof}
For $v \in E$ and $c \in \fld^\times$, replacing $v$ with $c v$
has the effect of multiplying $v \otimes \PhMap(v)$ by $c^2$, which does
not change the matroid $\DRig(M)$.

For the second assertion, let $g \in GL_r(\fld)$, and
suppose that we have an $\fld(\PhMap)$-linear dependence
  \begin{equation}
  \label{typical-dependence}
  \sum_{i=1}^n c_i v_i \otimes v^{(i)} = 0
  \end{equation}
in $\fld^r \otimes \fld(\PhMap)^d$.  The group $GL_r(\fld)$ acts
$\fld(\PhMap)$-linearly on $\fld^r \otimes \fld(\PhMap)^d$ by
$g(v \otimes w) = g(v) \otimes w$.  Applying $g$ to \eqref{typical-dependence} yields
  $$
  \sum_{i=1}^n c_i g(v_i) \otimes v^{(i)} = g(0) = 0.
  $$
Equivalently, 
  $$
  \sum_{i=1}^n c_i g(v_i) \otimes (\PhMap \circ g^{-1}) g(v_i) = 0.
  $$
The entries of the $d \times r$ matrix $\PhMap \circ g^{-1}$ are
algebraically independent transcendentals over $\fld$ (because $\PhMap$ was),
and the transcendental extensions $\fld(\PhMap)$ and $\fld(\PhMap \circ g^{-1})$
coincide because $g$ is invertible.
Hence the matroid represented by $\{g(v_1), \ldots, g(v_n)\}$
contains the same dependence \eqref{typical-dependence} as do $\{v_1,\ldots,v_n\}$.
Considering all such dependences and replacing $g$ with $g^{-1}$, one sees that
this matroid is combinatorially identical to $\DRig(M)$.
\end{proof}

\begin{qn} \label{representation-dependent}
Is $\DRig(M)$ a combinatorial invariant of $M$, or does it
depend on the choice of field $\fld$ and the particular representation $\{v_1,\dots,v_n\}$ 
of $M$ in $\fld^r$?
\end{qn}

\noindent
In the special case $d=2$, the Nesting Theorem implies that $\DRig(M)$
is indeed a combinatorial invariant.
While we have no reason to expect invariance in all cases,
we have not found a counterexample.  We have seen that $\DRig(M)$ is indeed
combinatorial when $M=\Umat{2}{3}$ or $\Umat{2}{4}$.  In what follows, we
describe a matroid with two projectively inequivalent representations
whose $d$-rigidity matroids coincide.

\begin{ex}\label{3-by-3}
Consider the following two sets of nine coplanar vectors in $\Rr^3$:
  \begin{align*}
  E  &= \{(1,0,0),(1,0,1),(1,0,2),(1,1,0),(1,1,1),(1,1,2),(1,2,0),(1,2,1),(1,2,2)\},\\
  E' &= \{(1,0,0),(1,0,1),(1,0,3),(1,2,0),(1,2,1),(1,2,3),(1,3,0),(1,4,1),(1,6,3)\}.
  \end{align*}

\unitlength=0.02in
\begin{center}
\begin{picture}(230,50)
\put( 5,5){\framebox(40,40){}}
\put(10,10){\makebox(0,0){$\bullet$}}
\put(10,25){\makebox(0,0){$\bullet$}}
\put(10,40){\makebox(0,0){$\bullet$}}
\put(25,10){\makebox(0,0){$\bullet$}}
\put(25,25){\makebox(0,0){$\bullet$}}
\put(25,40){\makebox(0,0){$\bullet$}}
\put(40,10){\makebox(0,0){$\bullet$}}
\put(40,25){\makebox(0,0){$\bullet$}}
\put(40,40){\makebox(0,0){$\bullet$}}
\put(25, 0){\makebox(0,0){$E$}}

\put(155,5){\framebox(70,40){}}
\put(160,10){\makebox(0,0){$\bullet$}}
\put(160,20){\makebox(0,0){$\bullet$}}
\put(160,40){\makebox(0,0){$\bullet$}}
\put(180,10){\makebox(0,0){$\bullet$}}
\put(180,20){\makebox(0,0){$\bullet$}}
\put(180,40){\makebox(0,0){$\bullet$}}
\put(190,10){\makebox(0,0){$\bullet$}}
\put(200,20){\makebox(0,0){$\bullet$}}
\put(220,40){\makebox(0,0){$\bullet$}}
\put(190, 0){\makebox(0,0){$E'$}}
\end{picture}
\end{center}

Let $M,M'$ be the matroids represented by $E,E'$ respectively.  These matroids
are combinatorially isomorphic, but $E$ and $E'$ are certainly projectively
inequivalent.  On the other hand, computations
using {\tt Mathematica} show that $\R^2(M)=\R^2(M')$ (= $\Umat{5}{9}$)
and that $\R^3(M)=\R^3(M')$ (the bases are the subsets of $E$
(resp.\ $E'$) of cardinality 6, except for the complements of
the eight affine lines.)
\end{ex}

We next discuss how $\R^d(M)$ stabilizes for large $d$.
Let $\omega:~ \fld(\PhMap_{1,1},\dots,\PhMap_{d+1,r}) 
\to \fld(\PhMap_{1,1},\dots,\PhMap_{d,r})$ be the map
sending $\PhMap_{d+1,j}$ to $0$ for every $j$.  Then $\omega$
takes linear dependences on
rows of $\R^{d+1}(M)$ to linear dependences on rows of $\R^d(M)$.
Therefore $\R^d(M) \subseteq \R^{d+1}(M)$.

Since there are only finitely many simplicial complexes on
a fixed finite ground set $E$, the tower
  $$
  M=\R^1(M) \subseteq \R^2(M) \subseteq \R^3(M) \subseteq \cdots
  $$
must eventually stabilize to some complex $\R^\infty(M)$.
We can say more precisely when this stabilization occurs.

\begin{prop}
\label{stabilization}
Let $M$ be a matroid represented by $E=\{v_1,\dots,v_n\}\subset\fld^r$,
where (without loss of generality) $M$ has rank $r$.
Then for every $d\geq r$,
  $$\R^d(M)=\R^r(M)=\R^\infty(M).$$
\end{prop}

\begin{proof}
Since $\R^d(M) \subseteq \R^{d+1}(M)$,
it suffices to prove that $\R^d(M) \subseteq \R^r(M)$ for $d \geq r$.
Let $\PhMap$ be an $r \times r$ matrix of transcendentals over $\fld$.
Suppose that we have a linear dependence of the form
\eqref{typical-dependence}.
Let $\psi$ be another $d \times r$ matrix of transcendentals,
so that $\fld(\PhMap) \hookrightarrow \fld(\PhMap,\psi)$ is a purely transcendental
extension.   Viewing the matrix $\psi$ as a $\fld(\PhMap,\psi)$-linear map, one can 
apply it to the second factor of $\fld^r \otimes \fld(\PhMap)^r$.  Applying this to
\eqref{typical-dependence} gives
  \begin{equation}
  \label{derived-dependence}
  \sum_{i=1}^n c_i v_i \otimes (\psi \circ \PhMap)(v_i) = 0,
  \end{equation}
which is an $\fld(\PhMap,\psi)$-linear dependence on the vectors
$\{v_i \otimes (\psi \circ \PhMap)(v_i)\}_{i=1,\ldots,n}$.

We claim that $\fld(\PhMap,\psi)$ is purely transcendental over $\fld(\psi\circ\PhMap)$.
To see this, first note that $\fld(\psi,\PhMap) = \fld(\psi\circ\PhMap,\PhMap^{-1})$.
That is, $\fld(\psi,\PhMap)$ can be obtained from $\fld(\psi\circ\PhMap)$
by adjoining $r^2$ elements, namely the entries of $\PhMap^{-1}$.  In particular,
the transcendence degree of $\fld(\psi,\PhMap)$ over $\fld(\psi\circ\PhMap)$
is at most $r^2$.
Similarly, the transcendence degree of $\fld(\psi\circ\PhMap)$ over $\fld$ is
at most $dr$.  But $\fld(\psi,\PhMap)$ clearly has transcendence degree $dr+r^2$
over $\fld$, and transcendence degree is additive in towers of field extensions
\cite[Thm.~VI.1.11]{Hungerford}, so both instances of ``at most'' may be replaced with
``exactly'', proving the claim.

By the existence of
the $\fld(\PhMap,\psi)$-linear dependence \eqref{derived-dependence},
we conclude that the vectors $\{v_i \otimes (\psi \circ \PhMap)(v_i)\}_{i=1,\ldots,n}$
must also be $\fld(\psi \circ \PhMap)$-linearly dependent.
Therefore $\R^d(M) \subseteq \R^r(M)$ as desired.
\end{proof}

When a matroid $M$ can be represented over different fields,
it is natural to ask how much $\R^d(M)$ can vary.  For instance,
if $M=M(G)$ is graphic, then the standard representation
\eqref{std-graphic} is valid over \emph{every} field $\fld$
and unique up to projective equivalence once the field is fixed,
as mentioned earlier.  For sufficiently large $d$, the $d$-rigidity matroid of $M(G)$ is
also independent of the choice of the field $\fld$, as we now explain.

\begin{prop}\label{graphic-boolean-prop}
Let $M=M(G)$ be the graphic matroid representing an $n$-vertex graph $G$ over an 
arbitrary field $\fld$, equipped with the standard representation
\eqref{std-graphic}.  Then $\R^n(M) = 2^E$ = $\R^\infty(M)$.
\end{prop}

\begin{proof}
Let $K_n$ be the complete graph on $n$ vertices. Since $R^n(M)$ is a 
row-selected submatrix of $R^n(M(K_n))$, it suffices to assume that $G=K_n$. 

To avoid overly cumbersome notation, we give the proof for $n=4$;
the argument for arbitrary $n$ should be clear from this case.
For $n=4$, the $6\x 12$ rigidity matrix $R^4(M(K_4))$ is as follows.
(Each nonzero entry is a binomial $\PhMap_{ij}-\PhMap_{ik}$, written on
two lines so that the matrix is not too wide for the page.)

\tiny
$$
\left[ \begin{array}{ccc|ccc|ccc|ccc}
\ontopof{\PhMap_{11}}{-\PhMap_{12}} & \ontopof{\PhMap_{21}}{-\PhMap_{22}} & \ontopof{\PhMap_{31}}{-\PhMap_{32}} &
\ontopof{\PhMap_{12}}{-\PhMap_{11}} & \ontopof{\PhMap_{22}}{-\PhMap_{21}} & \ontopof{\PhMap_{32}}{-\PhMap_{31}} & 0 & 0 & 0 & 0 & 0 & 0 \\&&&&&&&&&&&\\
\ontopof{\PhMap_{11}}{-\PhMap_{13}} & \ontopof{\PhMap_{21}}{-\PhMap_{23}} & \ontopof{\PhMap_{31}}{-\PhMap_{33}} & 0 & 0 & 0 &
\ontopof{\PhMap_{13}}{-\PhMap_{11}} & \ontopof{\PhMap_{23}}{-\PhMap_{21}} & \ontopof{\PhMap_{33}}{-\PhMap_{31}} & 0 & 0 & 0 \\&&&&&&&&&&& \\
0 & 0 & 0 & \ontopof{\PhMap_{12}}{-\PhMap_{13}} & \ontopof{\PhMap_{22}}{-\PhMap_{23}} & \ontopof{\PhMap_{32}}{-\PhMap_{33}} &
\ontopof{\PhMap_{13}}{-\PhMap_{12}} & \ontopof{\PhMap_{23}}{-\PhMap_{22}} & \ontopof{\PhMap_{33}}{-\PhMap_{32}} & 0 & 0 & 0 \\&&&&&&&&&&& \\
\ontopof{\PhMap_{11}}{-\PhMap_{14}} & \ontopof{\PhMap_{21}}{-\PhMap_{24}} & \ontopof{\PhMap_{31}}{-\PhMap_{34}} & 0 & 0 & 0 & 0 & 0 & 0 &
\ontopof{\PhMap_{14}}{-\PhMap_{11}} & \ontopof{\PhMap_{24}}{-\PhMap_{21}} & \ontopof{\PhMap_{34}}{-\PhMap_{31}} \\&&&&&&&&&&& \\
0 & 0 & 0 & \ontopof{\PhMap_{12}}{-\PhMap_{14}} & \ontopof{\PhMap_{22}}{-\PhMap_{24}} & \ontopof{\PhMap_{32}}{-\PhMap_{34}} & 0 & 0 & 0 &
\ontopof{\PhMap_{14}}{-\PhMap_{12}} & \ontopof{\PhMap_{24}}{-\PhMap_{22}} & \ontopof{\PhMap_{34}}{-\PhMap_{32}} \\&&&&&&&&&&& \\
0 & 0 & 0 & 0 & 0 & 0 & \ontopof{\PhMap_{13}}{-\PhMap_{14}} & \ontopof{\PhMap_{23}}{-\PhMap_{24}} & \ontopof{\PhMap_{33}}{-\PhMap_{34}} &
\ontopof{\PhMap_{14}}{-\PhMap_{13}} & \ontopof{\PhMap_{24}}{-\PhMap_{23}} & \ontopof{\PhMap_{34}}{-\PhMap_{33}}
\end{array}
\right]
$$

\normalsize
We must show that some $6\times 6$ minor of $R^4(M(K_4))$ is nonsingular.
Consider the submatrix $M'$ consisting of the last column in the second block, the last
two columns in the third block, and all three columns in the fourth block:


$$
\left[\begin{array}{c|cc|ccc}
\mathbold{\PhMap_{32}-\PhMap_{31}} & 0 & 0 & 0 & 0 & 0 \\ \hline
0 & \mathbold{\PhMap_{23}-\PhMap_{21}} & \mathbold{\PhMap_{33}-\PhMap_{31}} & 0 & 0 & 0 \\
\PhMap_{32}-\PhMap_{33} & \mathbold{\PhMap_{23}-\PhMap_{22}} & \mathbold{\PhMap_{33}-\PhMap_{32}} & 0 & 0 & 0 \\ \hline
0 & 0 & 0 & \mathbold{\PhMap_{14}-\PhMap_{11}} & \mathbold{\PhMap_{24}-\PhMap_{21}} & \mathbold{\PhMap_{34}-\PhMap_{31}} 
\\
\PhMap_{32}-\PhMap_{34} & 0 & 0 & \mathbold{\PhMap_{14}-\PhMap_{12}} & \mathbold{\PhMap_{24}-\PhMap_{22}} & 
\mathbold{\PhMap_{34}-\PhMap_{32}} \\
0 & \PhMap_{23}-\PhMap_{24} & \PhMap_{33}-\PhMap_{34} & \mathbold{\PhMap_{14}-\PhMap_{13}} & \mathbold{\PhMap_{24} - 
\PhMap_{23}} & \mathbold{\PhMap_{34}-\PhMap_{33}}
\end{array}
\right]
$$

Since $M'$ is block lower triangular, its determinant is the product
of the determinants of the blocks along the diagonal (indicated in boldface).
Each such determinant is a nonzero polynomial in the $\PhMap_{ij}$
over any field, because the coefficients of
$\PhMap_{31}$ in the first block,
$\PhMap_{21}\PhMap_{32}$ in the second block, and 
$\PhMap_{11}\PhMap_{22}\PhMap_{33}$ in the third block
are all $\pm 1$. Therefore $M'$ is
nonsingular over any field, as desired.
\end{proof}

This observation begs the question of whether
$\R^d(M(G))$ depends on the field \emph{before} $d$ reaches the stable range.
For an arbitrary representable matroid $M$, it is not true in
general that $\R^\infty(M)$ is Boolean.  We have already seen one example
for which this fails, namely $\Umat{2}{4}$.  Another
example is the well-known \defterm{Fano matroid} $F$,
represented over the two-element field $\fld_2$ by the seven
nonzero elements of $\fld_2^3$.  It is not hard to show that
$\LL^d(F)$ is Boolean for $d>\frac{7}{3}$.  On the other hand,
computation with {\tt Mathematica} indicates that
$\R^2(F)=\Umat{5}{7}$, but
$\R^d(F)=\Umat{6}{7}$ for all integers $d\geq 3$.



\section{Open problems}

The foregoing results raise many questions that we think are worthy of
further study; some of these have been mentioned earlier in the paper.
In this final section, we restate the open problems and add a few more.

\begin{problem}
Determine the singular locus of the $(k,d)$-photo space $X_{k,d}(M)$
(perhaps by calculating the dimension of its various tangent spaces,
as in Proposition~\ref{smooth-Boolean}).
\end{problem}

\begin{problem}
Give a direct combinatorial explanation for Corollary~\ref{curious-Tutte-symmetry},
presumably by identifying some  natural relationship between photos of $M$ and of $M^\perp$.
\end{problem}

\begin{problem}
Explain the ``scaling phenomenon'' of Remark~\ref{funky-equalities} geometrically.
\end{problem}

\begin{problem}
Determine whether or not the $d$-rigidity matroid
$\DRig(M)$ is a combinatorial invariant of $M$
(Question~\ref{representation-dependent}).
If not, determine which matroids have this property,
and to what extent $\DRig(M)$ depends on the
field $\fld$ over which $M$ is represented.  
In particular, is $\DRig(M)$ independent of $\fld$ in the case that $M$
is a graphic matroid?
\end{problem}

Crapo gave an elegant characterization
\cite[Theorem 8.2.2]{Whiteley} of $\HH^d(M)$ when $M$ is graphic.
A basis of $\HH^d(M)$ is a (multi-)set of edges having
a \defterm{$(d+1){\bf T}d$-covering}, or
a decomposition into $d+1$ edge-disjoint trees, 
exactly $d$ incident with each vertex, with no $d$ nonempty subtrees
spanning the same subset of vertices.

\begin{problem}
Generalize Crapo's characterization of $\HH^d(M)$ to the case of a non-graphic matroid $M$.
\end{problem}

\noindent A vertex of a graph $G$ corresponds to a cocircuit of $M(G)$ whose deletion
leaves a connected matroid.  However, there is no analogous notion of
``vertex'' when $M$ is a non-graphic matroid (although the foregoing
may be helpful if $M$ is sufficiently connected).
Similarly, it is unclear how to generalize to non-graphic matroids
(and to higher dimensions)
other fundamentals of graphic rigidity theory; for instance,
Henneberg's construction of the bases for $\HH^2(M)=\R^2(M)=\LL^2(M)$
\cite[Theorem 2.2.3]{Whiteley}.

Our last open problem is similar in spirit to the results of
\cite{JLM1} and \cite{JLM3}, describing the algebraic and
combinatorial structure of the equations defining
the slope variety of a graph.  It is motivated also by the appearance of
the cross-ratio in Example~\ref{U24-example}.

\begin{problem}
Describe explicitly the defining equations
(in Pl\"ucker coordinates on $\Gr(k,\fld^d)^n$)
for $\ov{\pi X^{\0}_{k,d}(M)}$,
where $\pi$ is the projection map of \eqref{projection-map}.
\end{problem}

\section*{Acknowledgments}

A substantial portion of this work was completed in July 2004 at the Park 
City Mathematical Institute, sponsored by the Institute for Advanced 
Study.  The authors particularly thank Walter Whiteley for sharing with them some of 
his manuscripts on classical rigidity theory.  They also thank
Gil Kalai, Ezra Miller, Neil White and G\"unter Ziegler for useful discussions,
and an anonymous referee for a careful reading and numerous thoughtful suggestions.


\end{document}